\documentclass{article}
\usepackage{amsmath}
\usepackage{mathrsfs}
\usepackage{txfonts}
\usepackage{bm}
\usepackage{mathrsfs}
\usepackage{latexsym,amsfonts,amssymb}
\begin{document}
\setcounter{page}{1}
\newtheorem{theorem}{Theorem}[section]
\newtheorem{fthm}[theorem]{Fundamental Theorem}
\newtheorem{remark}[theorem]{Remark}
\newtheorem{lemma}[theorem]{Lemma}
\newtheorem{corollary}[theorem]{Corollary}
\newtheorem{example}[theorem]{Example}
\newtheorem{proposition}[theorem]{Proposition}

\renewcommand{\thefootnote}{\fnsymbol{footnote}}
\newcommand{\qed}{\hfill\Box\medskip}
\newcommand{\proof}{{\it Proof.\quad}}
\newcommand{\rmnum}[1]{\romannumeral #1}
\renewcommand{\abovewithdelims}[2]{%
\genfrac{[}{]}{0pt}{}{#1}{#2}}

\renewcommand{\thefootnote}{\fnsymbol{footnote}}

\begin{center}{\LARGE\bf Graph homomorphisms on rectangular matrices \\ over
division rings I}\footnote{Project 11371072 supported by National Natural Science Foundation of China.}
\end{center}
\begin{center}
Li-Ping Huang\footnote{E-mail: \ lipingmath@163.com (L.P. Huang); zhaokangmath@126.com (K. Zhao).}, \ \ Kang Zhao
\\
{\small  School of Math. and Statis.,
 {Changsha University of Science and Technology, Changsha, 410004,  China}}
\end{center}

\begin{abstract}
Let $\mathbb{D}$ be a division ring, and let
${\mathbb{D}}^{m\times n}$  be the set of $m\times n$  matrices over $\mathbb{D}$.
Two matrices $A,B\in {\mathbb{D}}^{m\times n}$ are adjacent if  ${\rm rank}(A-B)=1$. By the adjacency, ${\mathbb{D}}^{m\times n}$ is a connected graph.
Suppose that $m,n,m',n'\geq2$ are integers and $\mathbb{D}'$ is a  division ring.
Using the weighted semi-affine map and algebraic method, we characterize  graph homomorphisms from ${\mathbb{D}}^{m\times n}$ to ${\mathbb{D}'}^{m'\times n'}$
(where $|\mathbb{D}|\geq 4$) under some weaker conditions.

\vspace{3mm}


\vspace{2mm}

\noindent{\bf 2010 AMS Classification}:  15A33, 51A10, 51D20, 05C60

\end{abstract}

\section{Introduction}

Throughout this paper, we assume that
 $\mathbb{D}, \mathbb{D}'$ are division rings,  ${\mathbb{D}}^*=\mathbb{D}\setminus \{0\}$, and $m, n, m',n'$ are positive integers.
  Denote by  $\mathbb{F}_q$ the finite field with $q$ elements where $q$ is a power of a prime.
Let $|X|$ be the cardinality  of a set $X$. Let ${\mathbb{D}}^n$ [resp. $^n{\mathbb{D}}$] denote the left [resp. right] vector space over
$\mathbb{D}$ whose elements are $n$-dimensional row [resp. column] vectors over  $\mathbb{D}$.
On the basic properties of matrices over a division ring, one may see literatures \cite{Thoms-Algebra,W3}.

Let ${\mathbb{D}}^{m\times n}$ and ${\mathbb{D}}^{m\times n}_r$ denote the sets of $m\times n$  matrices and  $m\times n$  matrices of rank $r$ over $\mathbb{D}$,
respectively. Denote the set of $n\times n$ invertible matrices over $\mathbb{D}$ by $GL_n(\mathbb{D})$.
Let $I_r$ ($I$ for short) be the $r\times r$ identity matrix,  $0_{m,n}$ the $m\times n$ zero matrix ($0$ for short) and $0_n=0_{n,n}$.
Let   $E^{m\times n}_{ij}$ ($E_{ij}$ for short) denote the $m\times n$ matrix whose $(i,j)$-entry is $1$ and all other entries are $0$'s.
Denote by  $^tA$ the transpose matrix of a matrix $A$. If $\sigma: \mathbb{D}\rightarrow \mathbb{D}'$ is a map and  $A=(a_{ij})\in {\mathbb{D}}^{m\times n}$,
we write  $A^{\sigma}=(a_{ij}^{\sigma})$ and $^tA^\sigma=\,^t(A^\sigma)$.
Let $(B_1,B_2, \ldots, B_k)$ denote an $m\times (n_1+n_2+\cdots +n_k)$  matrix, where $B_i$ is an $m\times n_i$ matrix.
 Denote by ${\rm diag}(A_1,\ldots, A_k)$ a block diagonal matrix where $A_i$ is an $m_i\times n_i$ matrix.

For $A, B\in  {\mathbb{D}}^{m\times n}$, $ad(A, B):={\rm rank}(A-B)$  is called
the {\em arithmetic distance} between $A$ and $B$. Clearly, $ad(A, B)\geq 0$, $ ad(A, B)=0 \Leftrightarrow A=B$; $ad(A, B)=ad(B, A)$;
$ ad(A, B)\leq ad(A, C)+ad(C, B)$. Two matrices $A$ and $B$ are said to be {\em adjacent},  denoted by $A\sim B$, if $ad(A, B)=1$.

A matrix in ${\mathbb{D}}^{m\times n}$ is also called a point.
Two distinct points  $A, B\in {\mathbb{D}}^{m\times n}$  are said to be of {\em distance} $r$, denoted by $d(A, B)=r$, if
 $r$ is the least positive integer for which there is a sequence
of $r+1$ points $X_0, X_1,\cdots,X_r\in {\mathbb{D}}^{m\times n}$ with $X_0=A$ and $X_r=B$, such that $X_0\sim X_1\sim \cdots\sim X_r$.
Define $d(A, A)=0$. It is well-known that (cf. \cite[Lemma 3.3]{LAA436Huang} or \cite{W3})
\begin{equation}\label{ad-d}
\mbox{$d(A, B)=ad(A,B)$, \ for all $A, B\in {\mathbb{D}}^{m\times n}$.}
\end{equation}

 In this paper, all graphs are undirected graphs without loops or multiple edges.
Let $\Gamma(\mathbb{D}^{m\times n})$ be the graph whose vertex set is ${\mathbb{D}}^{m\times n}$ and
two vertices $A,B\in {\mathbb{D}}^{m\times n}$ are adjacent if
$ad(A, B)=1$. The  $\Gamma(\mathbb{D}^{m\times n})$ is called  the graph on $m\times n$ matrices  over $\mathbb{D}$.
By (\ref{ad-d}), $\Gamma(\mathbb{D}^{m\times n})$ is a connected graph. Recall that every $A\in {\mathbb{D}}^{m\times n}$
can be written as $A=P{\rm diag}(I_r,0)Q$, where $P,Q$ are invertible matrices over $\mathbb{D}$.
It is easy to prove that $\Gamma(\mathbb{D}^{m\times n})$ is a {\em distance transitive graph} \cite{Godsil}.
When $\mathbb{D}=\mathbb{F}_q$ is a finite field, $\Gamma(\mathbb{F}_q^{m\times n})$ is also called a {\em bilinear forms graph}  \cite{Brouwera2}.

Let  $\varphi:  \mathbb{D}^{m\times n} \rightarrow  {\mathbb{D}'}^{m'\times n'}$ be a map.  The map  $\varphi$ is called a {\em graph homomorphism}
if $A\sim B$ implies that $\varphi(A)\sim\varphi(B)$.
  The map $\varphi$ is called a {\em graph isomorphism} if $\varphi$ is  bijective  and $A\sim B$ $\Leftrightarrow$ $\varphi(A)\sim \varphi(B)$.
 The  map $\varphi$ is called a {\em  distance preserving map}
if $d(A, B)=d(\varphi(A),\varphi(B))$ for all $A,B\in {\mathbb{D}}^{m\times n}$.
The map $\varphi$ is called  a {\em  distance $k$ preserving map} if $d(A, B)=k$ implies that $d(\varphi(A),\varphi(B))=k$ for some fixed $k$.
 In the geometry of matrices, a graph homomorphism [resp. {\em graph isomorphism}] is also called an {\em adjacency preserving map}
[resp. {\em adjacency preserving bijection in both directions}].
If $\varphi:  {\mathbb{D}}^{m\times n}\rightarrow  {\mathbb{D}'}^{m'\times n'}$ is a graph homomorphism, then (\ref{ad-d}) implies that
\begin{equation}\label{uy87mbmm}
\mbox{$d(\varphi(A), \varphi(B))\leq d(A, B)$, \ for all $A, B\in  {\mathbb{D}}^{m\times n}$.}
\end{equation}

 In the {\em geometry of matrices} \cite{W3},  all bijective maps of the form
$X\mapsto PXQ+A$ where $P,Q$ are invertible matrices over $\mathbb{D}$
and $A\in\mathbb{D}^{m\times n}$, form a transformation group on the space ${\mathbb{D}}^{m\times n}$. The fundamental problem in the geometry of  matrices
is to  characterize the transformation group  by as few geometric invariants as possible.
In 1951,  Hua \cite{Hua8,CGroups} showed that the invariant ``adjacency" is almost sufficient to characterize the transformation group,
 and  proved the  fundamental theorem of the geometry of rectangular matrices over division rings. Hua's theorem also characterized the graph isomorphisms on
 rectangular matrices, and his work was continued by many scholars
 (cf. \cite{Chooi}, \cite{Chooi2}, \cite{LAA436Huang}-\cite{HWEN}, \cite{Lim-2}-\cite{Wan5}).
 Thus, the basic problem of the geometry of  matrices is also  to study graph isomorphisms and graph homomorphisms on matrices.

In the algebraic graph theory, the study of graph homomorphism is an important subject (cf. \cite{Godsil,G.Hahn, Oxford,Knauer}).
Thus, it is of significance to determine  graph homomorphisms on matrices.
However,  graph homomorphisms from $\Gamma(\mathbb{D}^{m\times n})$ to $\Gamma({\mathbb{D}'}^{m'\times n'})$ still remain to be further solved.
Recently, literature \cite{Huang-Li} proved that every  graph endomorphism of  $\Gamma(\mathbb{F}_q^{m\times n})$ is either an automorphism
or a colouring.
\v{S}emrl \cite{optimal} characterized the graph homomorphism $\phi:  \mathbb{D}^{n\times n}\rightarrow  {\mathbb{D}}^{p\times q}$ which satisfies the condition
that $\phi(0)=0$ and there exists $A_0\in {\mathbb{D}}^{n\times n}$ such that ${\rm rank}(\phi(A_0))=n$ (where $p,q\geq n\geq 3$ and $|\mathbb{D}|\geq 4$).
Huang and \v{S}emrl \cite{Huang-Semrl-2014} further solved the remaining case by describing the graph homomorphism
$\phi:  \mathbb{D}^{2\times 2}\rightarrow  {\mathbb{D}}^{p\times q}$. When $\mathbb{D}$ is an EAS division ring, Pazzis and \v{S}emrl \cite{Pazzis}
showed  that every graph homomorphism $\phi: \mathbb{D}^{n\times m}\rightarrow  {\mathbb{D}}^{p\times q}$ is either  ``standard" or ``degenerate".
 These  interesting work encourage us to further characterize the graph homomorphisms  from $\mathbb{D}^{m\times n}$
 to ${\mathbb{D}'}^{m'\times n'}$.

 In this paper, by  the weighted semi-affine map and  algebraic method, we will characterize the non-degenerate graph homomorphisms
 from ${\mathbb{D}}^{m\times n}$ to ${\mathbb{D}'}^{m'\times n'}$ under some weaker conditions.

A nonempty subset $\cal S$ of ${\mathbb{D}}^{m\times n}$ is called an {\em adjacent set} if every pair of distinct points of $\cal S$ are adjacent.
 An adjacent set  in ${\mathbb{D}}^{m\times n}$ is a {\em maximal adjacent set} ({\em maximal set} for short),
 if  there is no adjacent set in ${\mathbb{D}}^{m\times n}$ which properly contains it as a subset.
 In the graph theory, a maximal  set is also called a {\em maximal clique} \cite{Brouwera2,Godsil}.
In  ${\mathbb{D}}^{m\times n}$, every adjacent set can be extended to a maximal set, and there are exactly two types of maximal sets.

Let
$${\cal M}_i=\left\{\sum_{j=1}^{n}x_jE_{ij}: x_j\in \mathbb{D}\right \}\subset {\mathbb{D}}^{m\times n}, \     i=1, \cdots, m;
$$
$${\cal N}_j=\left\{\sum_{i=1}^{m} y_iE_{ij}: y_i\in \mathbb{D}\right \}\subset {\mathbb{D}}^{m\times n},  \ j=1, \cdots,n; $$
 $$ \ \ \ \ \ {\cal M}'_i=\left\{\sum_{j=1}^{n'}x_jE_{ij}: x_j\in \mathbb{D}' \right \}\subset {\mathbb{D}'}^{m'\times n'}, \     i=1, \cdots, m';$$
$$\ \ \ \ {\cal N}'_j=\left\{\sum_{i=1}^{m'} y_iE_{ij}: y_i\in \mathbb{D}'\right \}\subset {\mathbb{D}'}^{m'\times n'},  \ j=1, \cdots,n'.$$
We say that ${\cal M}_1$ and ${\cal N}_1$ [resp. ${\cal M}'_1$ and ${\cal N}_1'$] are  {\em standard maximal sets} of type one and type two in
 ${\mathbb{D}}^{m\times n}$ [resp. ${\mathbb{D}'}^{m'\times n'}$], respectively. Clearly,  ${\cal M}_1$ [resp. ${\cal N}_1$] is a left [resp. right] linear space
 over $\mathbb{D}$.

\begin{lemma}\label{Rectangular-PID2-4}{\em(cf. \cite[Prpopsition 3.8]{W3})}
In ${\mathbb{D}}^{m\times n}$, all ${{\cal M}_i}$'s, ${\mathcal{N}_j}$'s,
$i=1, \cdots, m, j=1, \cdots, n$, are maximal sets. Moreover,  any maximal set  ${\cal M}$ is one of  the
following  forms.

{\em Type one.}  \ ${\cal M}=P{\cal M}_1+A$, where $P\in GL_m(\mathbb{D})$ and $A\in {\mathbb{D}}^{m\times n}$.

{\em Type two.} \ ${\cal M}={\cal N}_1Q+A$, where  $Q\in GL_n(\mathbb{D})$ and  $A\in {\mathbb{D}}^{m\times n}$.

\end{lemma}

\vspace{2mm}

  Let $\mathbb{D}_{\leq 1}^{m\times n}=\mathbb{D}_{1}^{m\times n}\cup \{0_{m,n}\}$.
For $A\in \mathbb{D}^{m\times n}$, $\mathbb{D}_{\leq 1}^{m\times n}+A$ is called  the {\em unit ball} with a central point $A$.
 If $f:  {\mathbb{D}}^{m\times n} \rightarrow  {\mathbb{D}'}^{m'\times n'}$ is a graph homomorphism, then
$f(\mathbb{D}_{\leq 1}^{m\times n}+A)\subseteq (\mathbb{D}_{\leq 1}^{m\times n}+f(A))$ for each $A\in\mathbb{D}^{m\times n}$.

  A graph homomorphism $f: {\mathbb{D}}^{m\times n}\rightarrow{\mathbb{D}'}^{m'\times n'}$ is called  {\em degenerate},
 if there exists a  matrix $A\in\mathbb{D}^{m\times n}_{\leq 1}$ and there are two  maximal sets $\mathcal{M}$ and $\mathcal{N}$ of
 different types in ${\mathbb{D}'}^{m'\times n'}$, such that $f\left(\mathbb{D}_{\leq 1}^{m\times n}+A\right)\subseteq \mathcal{M}\cup \mathcal{N}$
 with $f(A)\in\mathcal{M}\cap\mathcal{N}$.
 A graph homomorphism $f$ is called {\em non-degenerate} if it is not degenerate.
Note that the definition  is different from literatures \cite{Pazzis, optimal}.
 A graph homomorphism $f: {\mathbb{D}}^{m\times n}\rightarrow{\mathbb{D}'}^{m'\times n'}$ is called a  {\em colouring} if $f(\mathbb{D}^{m\times n})$ is an adjacent set.
Every colouring is degenerate.

A non-degenerate  graph homomorphism may not preserve distance $2$ (see below Example \ref{exatfds564776s}).
If $f$ is a degenerate [resp. non-degenerate] graph homomorphism, then  the map
$X\longmapsto P_1f(Q_1XQ_2)P_2$
is also a degenerate [resp. non-degenerate] graph homomorphism, where $Q_1,Q_2$ [resp. $P_1,P_2$] are invertible matrices over $\mathbb{D}$ [resp. $\mathbb{D}'$].

This  paper is organized as follows.
In Section 2, we  statement our main result and some examples on the graph homomorphism. We also give a sufficient and necessary condition on existence of the graph homomorphism.
In Section 3, we will discuss  maximal sets on ${\mathbb{D}}^{m\times n}$ and weighted semi-affine maps.
In Section 4, we will discuss Grassmann spaces and related lemmas.
In Section 5, we will prove our main result.

\section{Statement of Main Result and Examples}

 Dimension is a basic concept of geometry. We define the dimension of an adjacent set as follows.

Let $\mathcal{S}$ be an adjacent set in ${\mathbb{D}}^{m\times n}$ with $0\in \mathcal{S}$ and $|\mathcal{S}|\geq 2$. Then there exists a maximal set $\mathcal{M}$
 containing $\mathcal{S}$.
By Lemma \ref{Rectangular-PID2-4}, $\mathcal{M}=P\mathcal{M}_1$ or $\mathcal{M}=\mathcal{N}_1Q$,  where $P$ and $Q$ are invertible matrices. Hence
$P^{-1}\mathcal{S}\subseteq\mathcal{M}_1$ or $\mathcal{S}Q^{-1}\subseteq\mathcal{N}_1$. When $P^{-1}\mathcal{S}\subseteq\mathcal{M}_1$
[resp. $\mathcal{S}Q^{-1}\subseteq\mathcal{N}_1$],
the {\em dimension} of $\mathcal{S}$,  denoted by ${\rm dim}(\mathcal{S})$,  is the
number of matrices in a maximal left [resp. right] linear independent subset  of $P^{-1}\mathcal{S}$ [resp. $\mathcal{S}Q^{-1}$].
The ${\rm dim}(\mathcal{S})$ is uniquely determined by $\mathcal{S}$ and ${\rm dim}(\mathcal{S})\leq {\rm max}\{m,n\}$. Moreover,
\begin{equation}\label{vxc53gdnjl66}
\mbox{${\rm dim}(\mathcal{S})={\rm dim}(P_1\mathcal{S}Q_1)$, \ where $P_1\in GL_m(\mathbb{D})$ and $Q_1\in GL_n(\mathbb{D})$.}
\end{equation}

A division ring $\mathbb{D}$ is said to be an {\em EAS division ring} if every nonzero endomorphism of $\mathbb{D}$ is automatically surjective.
The following division rings are EAS division rings: the field of real numbers, the field of rational numbers and the finite fields.


Now, we statement our main result  as follows.

\begin{theorem}\label{MainTheorem0}   Let $\mathbb{D}, \mathbb{D}'$ be division rings with $|\mathbb{D}|\geq 4$,  and let $m,n,m', n'\geq 2$ be
integers. Suppose that $\varphi: {\mathbb{D}}^{m\times n}\rightarrow  {\mathbb{D}'}^{m'\times n'}$ is a non-degenerate graph homomorphism  with  $\varphi(0)=0$.
Assume further that ${\rm dim}(\varphi(\mathcal{M}_1))=n$ and  ${\rm dim}(\varphi(\mathcal{N}_1))=m$.
Then either there exist  matrices  $P\in GL_{m'}(\mathbb{D}')$ and $Q\in GL_{n'}(\mathbb{D}')$,
a nonzero ring homomorphism $\tau: \mathbb{D} \rightarrow\mathbb{D}'$, and a matrix $L\in {\mathbb{D}'}^{n\times m}$ with the property that $I_m+X^\tau L\in GL_m(\mathbb{D}')$
for every $X\in {\mathbb{D}}^{m\times n}$, such that $m'\geq m$, $n'\geq n$ and
\begin{equation}\label{tt3654665cnc00}
 \varphi(X)=P\left(\begin{array}{cc}
 (I_m+X^\tau L)^{-1}X^\tau  & 0 \\
  0 & 0 \\
   \end{array}
  \right)Q, \ \, X\in {\mathbb{D}}^{m\times n};
\end{equation}
or there exist matrices  $P\in GL_{m'}(\mathbb{D}')$ and $Q\in GL_{n'}(\mathbb{D}')$,
a nonzero ring anti-homomorphism $\sigma: \mathbb{D}\rightarrow\mathbb{D}'$, and a matrix $L\in {\mathbb{D}'}^{m\times n}$
with the property that $I_m+L\,^tX^\sigma\in GL_m(\mathbb{D}')$ for every $X\in {\mathbb{D}}^{m\times n}$, such that $m'\geq n$, $n'\geq m$ and
\begin{equation}\label{ttd2CX3mmbb6600}
\varphi(X)=P\left(\begin{array}{cc}
{^tX^\sigma}(I_m+L\,^tX^\sigma)^{-1}  & 0 \\
0 & 0 \\
\end{array}
\right)Q, \ \,  X\in {\mathbb{D}}^{m\times n}.
\end{equation}
In particular, if $\tau$ $[$resp. $\sigma$$]$ is surjective, then $L=0$ and  $\tau$ is a ring isomorphism
$[$resp. $\sigma$ is a ring anti-isomorphism$]$.
In {\rm(\ref{tt3654665cnc00})} or {\rm(\ref{ttd2CX3mmbb6600})}, some zero elements of right matrix may be absent.
\end{theorem}

\begin{corollary}\label{MainTheoremcor0}{\em(cf. \cite{Pazzis, optimal})} \  Let $\mathbb{D}$ be an EAS division ring with $|\mathbb{D}|\geq 4$,  and let $m,n,m', n'\geq 2$ be
integers. Suppose that $\varphi: {\mathbb{D}}^{m\times n}\rightarrow  {\mathbb{D}}^{m'\times n'}$ is a non-degenerate graph homomorphism with $\varphi(0)=0$.
Then either there exist  matrices  $P\in GL_{m'}(\mathbb{D})$ and $Q\in GL_{n'}(\mathbb{D})$,
and a ring automorphism $\tau$ of $\mathbb{D}$,  such that $m'\geq m$, $n'\geq n$ and
\begin{equation}\label{t3654665cnc00}
 \varphi(X)=P\left(\begin{array}{cc}
X^\tau  & 0 \\
  0 & 0 \\
   \end{array}
  \right)Q, \ \, X\in {\mathbb{D}}^{m\times n};
\end{equation}
or there exist matrices  $P\in GL_{m'}(\mathbb{D})$ and $Q\in GL_{n'}(\mathbb{D})$,
and a ring anti-automorphism $\sigma$ of $\mathbb{D}$,  such that $m'\geq n$, $n'\geq m$ and
\begin{equation}\label{d2CX3mmbb6600}
\varphi(X)=P\left(\begin{array}{cc}
{^tX^\sigma} & 0 \\
0 & 0 \\
\end{array}
\right)Q, \ \,  X\in {\mathbb{D}}^{m\times n}.
\end{equation}
\end{corollary}

We give the following  examples on graph homomorphisms.

\begin{example}\label{exatfds564776s}   There exists a non-degenerate graph homomorphism which does not preserve distance $2$. In fact,
 let $\mathbb{D}, \mathbb{D}'$ be division rings such that there exists $0\neq \xi\in\mathbb{D}'\setminus \mathbb{D}$.
Write $\scriptsize P_1=\left(
           \begin{array}{ccc}
             1 & 0& \xi \\
             0 & 1 & \xi \\
           \end{array}
         \right)$.
Define the map $\varphi: \mathbb{D}^{3\times n}\rightarrow {\mathbb{D}'}^{3\times n}$ $(n\geq 2)$ by $\scriptsize\varphi(X)=\left(
             \begin{array}{c}
               P_1X \\
               0 \\
             \end{array}
           \right)$, $X\in \mathbb{D}^{3\times n}$. Then
\begin{equation}\label{bvcrd25gjrfs87}
\varphi\left(
         \begin{array}{c}
           x \\
           y \\
           z \\
         \end{array}
       \right)
=\left(
             \begin{array}{c}
               x+\xi z \\
               y+\xi z \\
               0\\
             \end{array}
           \right), \ x,y,z\in \mathbb{D}^{n}.
\end{equation}
The map $\varphi$ is an additive map and a graph homomorphism.
Since $\scriptsize\varphi\left(
         \begin{array}{c}
           x \\
           x \\
           z \\
         \end{array}
       \right)
=\left(
             \begin{array}{c}
               x+\xi z \\
               x+\xi z \\
               0\\
             \end{array}
           \right)$, $\varphi$ does not preserve distance $2$.
We assert that $\varphi$ is non-degenerate. In fact,
 since $\varphi$ is additive, $\varphi\left(\mathbb{D}_{\leq 1}^{3\times n}+A\right)
=\varphi\left(\mathbb{D}_{\leq 1}^{3\times n}\right)+\varphi(A)$ for any $A\in\mathbb{D}^{3\times n}$. We have that $\varphi(0)=0$, $\varphi(\mathcal{M}_i)\subseteq \mathcal{M}_i'$, $i=1,2$,
and $\varphi(\mathcal{N}_j)\subseteq \mathcal{N}_j'$, $j=1,\ldots, n$. Thus, there are no two maximal sets $\mathcal{M}, \mathcal{N}$ of
 different types in ${\mathbb{D}'}^{3\times n}$, such that $0\in\mathcal{M}\cap \mathcal{N}$ and
 $\varphi\left(\mathbb{D}_{\leq 1}^{3\times n}\right)\subseteq \mathcal{M}\cup \mathcal{N}$.
 It follows that $\varphi\left(\mathbb{D}_{\leq 1}^{3\times n}+A\right)$ is not contained in a union of two maximal sets of
  different types containing $\varphi(A)$.
 Therefore, the  $\varphi$ is a non-degenerate graph homomorphism which  does not preserve distance $2$.
\end{example}

\begin{example}\label{exa-infinite}  Let $m,n,m', n'\geq 2$ be integers, and let
 $\mathbb{D}$ and $\mathbb{D}'$ be division rings such that $|\mathbb{D}|\leq|\mathbb{D}'|$ and $|\mathbb{D}|$ is infinite.
Then there is a colouring from $\mathbb{D}^{m\times n}$ to ${\mathbb{D}'}^{m'\times n'}$.
In fact, by the cardinal numbers theory (cf. \cite[Introduction]{Thoms-Algebra}),
we have $|\mathbb{D}|=|\mathbb{D}^{m\times n}|\leq|\mathbb{D}'|$, and hence there is an injective map
$f: \mathbb{D}^{m\times n}\rightarrow {\mathbb{D}'}^{1\times 1}$ such that $f$ is a colouring.
 Since there exists an embedding map (or inclusion mapping) ${\mathbb{D}'}^{1\times 1}\rightarrow {\mathbb{D}'}^{m'\times n'}$ which is a graph homomorphism,
  there is a colouring  form $\mathbb{D}^{m\times n}$ to ${\mathbb{D}'}^{m'\times n'}$.
\end{example}

However, we cannot write  algebraic formulas of the graph homomorphisms in the case of Example \ref{exa-infinite}.
In order to write the algebraic formula of a graph homomorphism we need some additional conditions.

\begin{example}\label{wett8345}{\em(cf. \cite{Hua8,CGroups,optimal})} \ Let $\mathbb{D}, \mathbb{D}'$ be division rings and let $m, n, m',n'\geq 2$ be integers.
Suppose that $\varphi: {\mathbb{D}}^{m\times n}\rightarrow  {\mathbb{D}'}^{m'\times n'}$ is a distance $k$ preserving map where $1\leq k\leq {\rm mim}\{m,n\}$, and
there is a fixed matrix $L\in {\mathbb{D}'}^{n'\times m'}$ with the property that $I_{m'}+\varphi(X) L\in GL_{m'}(\mathbb{D}')$ for every $X\in {\mathbb{D}}^{m\times n}$.
Put
\begin{equation}\label{bc5358gitlfddg}
\theta(X)=(I_{m'}+\varphi(X) L)^{-1}\varphi(X), \ X\in {\mathbb{D}}^{m\times n}.
\end{equation}
Then the map $\theta:  {\mathbb{D}}^{m\times n}\rightarrow  {\mathbb{D}'}^{m'\times n'}$ is a  distance $k$ preserving map.
\end{example}
\proof
For any $X,Y\in {\mathbb{D}}^{m\times n}$ with $d(X,Y)=k$, we have  $k={\rm rank}(X-Y)=d(\varphi(X),\varphi(X))={\rm rank}(\varphi(X)-\varphi(Y))$. Thus
\begin{eqnarray*}
  && m'+{\rm rank}(\varphi(X)-\varphi(Y))= {\rm rank}\left[\left(
\begin{array}{cc}
I_{m'} & \varphi(X) \\
I_{m'} & \varphi(Y)\\
\end{array}
\right)\left(
\begin{array}{cc}
I_{m'} & 0 \\
L & I_{n'} \\
\end{array}
\right) \right]\\
&=&{\rm rank}\left(
\begin{array}{cc}
I_{m'}+\varphi(X) L & \varphi(X) \\
I_{m'}+\varphi(Y) L & \varphi(Y)\\
\end{array}
\right)
= {\rm rank}\left(
\begin{array}{cc}
I_{m'} & \theta(X) \\
I_{m'} &  \theta(Y)\\
\end{array}
\right)=m'+d(\theta(X),\theta(Y)).
\end{eqnarray*}
Hence $d(X,Y)=d(\theta(X),\theta(Y))=k$. It follows that $\theta$ is a distance $k$ preserving map.
$\qed$

Similarly, we have

\begin{example}\label{wettwq345}{\rm (cf. \cite{Hua8,CGroups,optimal})} \ Let $\mathbb{D}, \mathbb{D}'$ be division rings and let $m, n,m', n'\geq 2$ be integers.
Suppose that $\varphi: {\mathbb{D}}^{m\times n}\rightarrow  {\mathbb{D}'}^{m'\times n'}$ is a distance $k$ preserving map where $1\leq k\leq {\rm mim}\{m,n\}$, and
there is a fixed matrix $L\in {\mathbb{D}'}^{n'\times m'}$ with the property that $I_{n'}+L\varphi(X)\in GL_{n'}(\mathbb{D}')$ for every $X\in {\mathbb{D}}^{m\times n}$.
Put
\begin{equation}\label{bc5358gitlfddg}
\theta(X)=\varphi(X)(I_{n'}+L\varphi(X))^{-1}, \ X\in {\mathbb{D}}^{m\times n}.
\end{equation}
Then the map $\theta:  {\mathbb{D}}^{m\times n}\rightarrow  {\mathbb{D}'}^{m'\times n'}$ is a  distance $k$ preserving map.
\end{example}

\begin{example}\label{Non-existence1}
Assume that  $m,n,m', n'\geq 1$ are integers, $|\mathbb{D}|>|\mathbb{D}'|$ and $|\mathbb{D}|$ is infinite.
Then there is no a graph homomorphism from $\mathbb{D}^{m\times n}$ to ${\mathbb{D}'}^{m'\times n'}$,
In fact, if there is a graph homomorphism $\varphi$ from $\mathbb{D}^{m\times n}$ to ${\mathbb{D}'}^{m'\times n'}$, then $\varphi(\mathcal{M}_1)\subseteq \mathcal{M}'$
where $\mathcal{M}'$ is a maximal set in ${\mathbb{D}'}^{m'\times n'}$. Since the restriction map
$f\mid_{\mathcal{M}_1}: \mathcal{M}_1\rightarrow \mathcal{M}'$ is injective, $\mathbb{D}'$ is also infinite.
By the cardinal numbers theory \cite{Thoms-Algebra},
we have $|\mathbb{D}'|=|\mathcal{M}'|\geq |\mathbb{D}|=|\mathcal{M}_1|$, a contradiction.
\end{example}

Thus, when $|\mathbb{D}|$ is infinite,  Example \ref{exa-infinite} implies that there is a graph homomorphism from $\mathbb{D}^{m\times n}$ to ${\mathbb{D}'}^{m'\times n'}$
if and only if $|\mathbb{D}|\leq|\mathbb{D}'|$.

For a finite graph $G$,
recall that a {\em maximum clique} of $G$ is a clique (i.e. adjacent subset) of $G$ which has maximum cardinality.
For convenience, we think that a maximum  (maximal) clique and its vertex set are equal. The {\em clique number} of $G$, denoted by  $\omega(G)$,  is
 the number of vertices in a maximum clique of $G$. Let $K_r$ be the {\em complete graph} on $r$ vertices.
Recall that an {\em $r$-colouring} of $G$ is a homomorphism from $G$ to $K_r$.
The {\em chromatic number} $\chi(G)$ of $G$ is  the least value $k$ for which $G$ can be  $k$-coloured (cf. \cite{Godsil}-\cite{Oxford}).

  The following is a necessary and sufficient condition on existence of graph homomorphism.

\begin{theorem}\label{Non-existence2}  Let $\mathbb{D}, \mathbb{D}'$ be division rings, and let $m,n,m', n'\geq 1$ be integers.
Put $s={\rm max}\{m,n\}$ and  $k={\rm max}\{m',n'\}$.
 Then, there exists a graph homomorphism from $\mathbb{D}^{m\times n}$ to ${\mathbb{D}'}^{m'\times n'}$
if and only if $|\mathbb{D}|^s\leq|\mathbb{D}'|^k$.
\end{theorem}
\proof
Let $\Gamma=\Gamma(\mathbb{D}^{m\times n})$ and $\Gamma'=\Gamma(\mathbb{{D}'}^{m'\times n'})$.
Without loss of generality, we assume that $n={\rm max}\{m,n\}$ and  $n'={\rm max}\{m',n'\}$.

Suppose $f: \mathbb{D}^{m\times n}\rightarrow {\mathbb{D}'}^{m'\times n'}$ is a graph homomorphism.
Then $\varphi(\mathcal{M}_1)\subseteq \mathcal{M}'$ where $\mathcal{M}'$ is a maximal set in ${\mathbb{D}'}^{m'\times n'}$. By Lemma \ref{Rectangular-PID2-4},
$\mathcal{M}'=P\mathcal{M}'_1+A$ or $\mathcal{M}'=\mathcal{N}'_1Q+A$, where $P\in GL_{m'}(\mathbb{D}')$,  $Q\in GL_{n'}(\mathbb{D}')$ and $A\in {\mathbb{D}'}^{m'\times n'}$ are fixed.
Clearly, $|{\mathcal{M}'}|\leq|\mathbb{D}'|^{n'}$.
Note that the restriction map $f\mid_{\mathcal{M}_1}: \mathcal{M}_1\rightarrow \mathcal{M}'$ is injective.
 By the cardinal numbers theory \cite{Thoms-Algebra}, we have $|{\mathcal{M}_1}|=|\mathbb{D}|^n\leq|{\mathcal{M}'}| \leq|\mathbb{D}'|^{n'}$.
Thus $|\mathbb{D}|^n\leq|\mathbb{D}'|^{n'}$.

Now, assume that $|\mathbb{D}|^n\leq|\mathbb{D}'|^{n'}$. We show that there exists a graph homomorphism  from $\mathbb{D}^{m\times n}$ to ${\mathbb{D}'}^{m'\times n'}$ as follows.

{\bf Case 1}. \  Both $|\mathbb{D}|$ and $|\mathbb{D}'|$ are finite.
By \cite[Theorem 2.9]{Huang-Li}, we have that $\chi(\Gamma)=|\mathbb{D}|^n=\omega(\Gamma)=|\mathcal{M}_1|$
and $\chi(\Gamma')=|\mathbb{D}'|^{n'}=\omega(\Gamma')=|\mathcal{M}'_1|$.
Let $\chi_1=|\mathbb{D}|^n$ and $\chi_2=|\mathbb{D}'|^{n'}$. Since $K_{\chi_2}$ and $\mathcal{M}'_1$ are isomorphic graphs,
we have graph homomorphisms $K_{\chi_2}\rightarrow \mathcal{M}'_1\rightarrow \Gamma'$, where the graph homomorphism
$\mathcal{M}'_1\rightarrow \Gamma'$ is an embedding map (or inclusion mapping).
By the algebraic graph theory (cf. \cite[p.104]{Godsil}), we have $\chi_1\leq \chi_2$. Thus  there are  graph homomorphisms
$\Gamma\rightarrow K_{\chi_1}\rightarrow K_{\chi_2}\rightarrow \mathcal{M}'_1\rightarrow \Gamma'$,
 where the graph homomorphism $K_{\chi_1}\rightarrow K_{\chi_2}$ is an embedding map.
Therefore,
there exists a graph homomorphism from $\mathbb{D}^{m\times n}$ to ${\mathbb{D}'}^{m'\times n'}$.

{\bf Case 2}. \  $|\mathbb{D}|$ is infinite. Then $|\mathbb{D}|^n=|\mathbb{D}|$. Since $|\mathbb{D}|^n\leq|\mathbb{D}'|^{n'}$,
 $|\mathbb{D}'|$ is infinite and  $|\mathbb{D}'|^{n'}=|\mathbb{D}'|$. Thus $|\mathbb{D}|\leq |\mathbb{D}'|$. By Example \ref{exa-infinite},
there is a graph homomorphism from $\mathbb{D}^{m\times n}$ to ${\mathbb{D}'}^{m'\times n'}$.

{\bf Case 3}. \  $|\mathbb{D}|$ is finite and $|\mathbb{D}'|$ is infinite. Then   $|\mathbb{D}'|=|\mathbb{D}'|^{n'}>|\mathbb{D}|^{mn}$.
Clearly, there is  a graph homomorphism from $\mathbb{D}^{m\times n}$ to ${\mathbb{D}'}^{1\times 1}$. On the other hand, there exists an embedding map
${\mathbb{D}'}^{1\times 1}\rightarrow {\mathbb{D}'}^{m'\times n'}$ which is a graph homomorphism.  Hence there is a graph homomorphism
from $\mathbb{D}^{m\times n}$ to ${\mathbb{D}'}^{m'\times n'}$.
$\qed$

\section{Maximal Sets and Weighted Semi-affine Maps}

In this section, we will discuss maximal sets on ${\mathbb{D}}^{m\times n}$ and their affine geometries.

\begin{lemma}\label{Matrix-PID5-4bb}{\em (cf. \cite[Lemma 3.2]{LAA436Huang})}  Let  $m,n,r,s$ be integers with $1\leq r,s<{\rm min}\{m, n\}$.
Assume that $\alpha=\{i_1, \ldots, i_r\}$, $\beta=\{j_1, \ldots, j_s\}$, where $1\leq i_1<\cdots <i_r\leq m$, $1\leq j_1<\cdots <j_s\leq n$.
Let $A=(a_{ij})\in \mathbb{D}^{m \times n}$,
$B_i=\sum_{t=1}^r\sum_{k=1}^sb^{(i)}_{i_tj_k}E_{i_tj_k}\in \mathbb{D}^{m \times n}$ (where $b^{(i)}_{i_tj_k}\in \mathbb{D}$), $i=1, 2$, and $B_1\neq B_2$.
If $A\sim B_i$, $i=1, 2$, then either $a_{ij}=0$ for all $i\notin \alpha$, or  $a_{ij}=0$ for all $j\notin \beta$.
\end{lemma}

Using Lemmas \ref{Matrix-PID5-4bb} and \ref{Rectangular-PID2-4}, we can prove the following results.

\begin{corollary}\label{Rectangular-PID2-7}{\em (cf. \cite[Corollary 3.10]{W3})} \
 Let $A$ and $B$ be two  adjacent points in ${\mathbb{D}}^{m\times n}$. Then there are exactly two maximal sets $\mathcal{M}$
 and $\mathcal{M}'$ containing $A$ and $B$. Moreover, $\mathcal{M}$ and $\mathcal{M}'$ are of different types.
\end{corollary}

\begin{corollary}\label{Rectangular-1-13}{\em (cf. \cite{LAA436Huang,Huang-book,W3})} \
If $\mathcal{M}$ and $\mathcal{N}$ are two distinct maximal sets of the same type $[$resp.  different types$]$ in ${\mathbb{D}}^{m\times n}$ with $\mathcal{M}\cap\mathcal{N}\neq \emptyset$,
 then $|\mathcal{M}\cap\mathcal{N}|=1$ $[$resp. $|\mathcal{M}\cap \mathcal{N}|=|\mathbb{D}|\geq 2$$]$.
\end{corollary}

\begin{lemma}\label{maximalset022} Suppose  $\mathcal{M}$ and $\mathcal{N}$ are two distinct maximal sets in  ${\mathbb{D}}^{m\times n}$
with $\mathcal{M}\cap\mathcal{N}\neq \emptyset$. Then:
\begin{itemize}
\item[{\em (i)}]  if $\mathcal{M}$ and $\mathcal{N}$ are of different types, then for any $A\in\mathcal{M}\cap\mathcal{N}$,
there are $P\in GL_m(\mathbb{D})$ and $Q\in GL_n(\mathbb{D})$
such that ${\cal M}=P{\cal M}_1Q+A=P{\cal M}_1+A$ and  ${\cal N}=P{\cal N}_1Q+A={\cal N}_1Q+A$;

\item[{\em (ii)}]  if both $\mathcal{M}$ and $\mathcal{N}$ are of type one $[$resp. type two$]$, then there exists an invertible  matrix $P$ $[$resp. $Q$$]$
such that ${\cal M}=P{\cal M}_1+A$ and ${\cal N}=P{\cal M}_2+A$
$[$resp. ${\cal M}={\cal N}_1Q+A$ and ${\cal N}={\cal N}_2Q+A$$]$, where $\mathcal{M}\cap\mathcal{N} = \{A \}$.
\end{itemize}
\end{lemma}
\proof
(i). By Lemma \ref{Rectangular-PID2-4}, it is easy to see that (i) holds.

(ii).  Let $\mathcal{M}$ and $\mathcal{N}$ be of type one with $\mathcal{M}\cap\mathcal{N}\neq \emptyset$.
By Corollary \ref{Rectangular-1-13}, $\mathcal{M}\cap\mathcal{N} = \{A \}$. By Lemma \ref{Rectangular-PID2-4}, there are $P_1, P_2\in GL_m(\mathbb{D})$
such that $\mathcal{M}=P_1\mathcal{M}_1+A$ and $\mathcal{N}=P_2\mathcal{M}_1+A$. Write $Q_1=P_1^{-1}P_2$. Then $\mathcal{N}=P_1Q_1\mathcal{M}_1+A$.
Since $\mathcal{M}\neq\mathcal{N}$, $\mathcal{M}_1\neq Q_1\mathcal{M}_1$. Clearly,
$$\mbox{$\mathcal{M}_1=\left\{\left(
                    \begin{array}{c}
                      x \\
                      0 \\
                    \end{array}
                  \right): x\in \mathbb{D}^n\right\}$, \ \ \
$Q_1\mathcal{M}_1=\left\{\left(
                    \begin{array}{c}
                      q_1x \\
                      \vdots \\
                      q_mx \\
                    \end{array}
                  \right): x\in \mathbb{D}^n\right\}$,}$$
where $^t(q_1,\ldots,q_m)\in \,^m\mathbb{D}$ and $^t(q_2,\ldots,q_m)\neq 0$.
Thus, there exists a $Q_2\in GL_m(\mathbb{D})$ such that   $Q_2Q_1\mathcal{M}_1=\mathcal{M}_2$ and $Q_2\mathcal{M}_1=\mathcal{M}_1$.
Put $P=P_1Q_2^{-1}$. We obtain that ${\cal M}=P{\cal M}_1+A$ and ${\cal N}=P{\cal M}_2+A$.

Similarly, if both $\mathcal{M}$ and $\mathcal{N}$ are  of type two with $\mathcal{M}\cap\mathcal{N} = \{A \}$, then there exists an invertible  matrix $Q$
such that ${\cal M}={\cal N}_1Q+A$ and ${\cal N}={\cal N}_2Q+A$.
$\qed$

Now, we simply introduce  the affine geometry.  Let $V$ be
 an $r$-dimensional  left vector subspace  of ${\mathbb{D}}^n$ and $a\in\mathbb{D}^n$. Then $V+a$ is called an $r$-dimensional {\em left affine flat} ({\em flat} for short)
over $\mathbb{D}$.  When $r\geq 2$, the set of all flats in $V+a$ is called the {\em left  affine geometry} on $V+a$, which is denoted by $AG(V+a)$.  The {\em dimension} of $AG(V+a)$ is $r$,
 denoted by ${\rm dim}(AG(V+a))=r$. The flats of dimensions $0, 1, 2$  are called {\em points, lines, planes} in $AG(V+a)$.
Similarly, we can define the {\em right  affine geometry} on a right affine flat over $\mathbb{D}$.
There is a general  axiomatic definition on affine geometry (cf. \cite{Bennett, ModernProjective}).

In a left or right  affine geometry, the {\em join} of flats $M_1$ and $M_2$,  denoted by $M_1\vee M_2$, is the minimum dimensional  flat containing  $M_1$ and $M_2$.
 The join $M_1\vee M_2$ is also the intersection of all flats containing $M_1$ and $M_2$ in the affine geometry.
 Thus, in an affine geometry, a line $\ell$ can be written as  $\ell =b\vee c$, where $b,c\in V$ are two distinct points.

Let ${\cal M}=P{\cal M}_1+A$ be a maximal set of type one in $\mathbb{D}^{m\times n}$, where $P\in GL_m(\mathbb{D})$ and $A\in\mathbb{D}^{m\times n}$. Then we have
a left affine geometry $AG(P{\cal M}_1+A)$ such that $AG(P{\cal M}_1+A)$ and $AG({\mathbb{D}}^n)$ are affine isomorphic.
Similarly,  we have a right affine geometry $AG({\cal N}_1Q+A)$ such that $AG({\cal N}_1Q+A)$ and $AG(^m{\mathbb{D}})$ are affine isomorphic,
where $Q\in GL_n(\mathbb{D})$.
In $AG(P{\cal M}_1+A)$, the parametric equation of a line $\ell$ is
$$\mbox{$\ell=P\left(
\begin{array}{cc}
\mathbb{D}\alpha+\beta \\
0\\
\end{array}\right)+A$, \ where $\alpha, \beta\in \mathbb{D}^n$ with $\alpha\neq 0$.}
$$

\begin{lemma}{\em (cf. \cite[p.95]{W3})}\label{Rectangular-PID2-11} \ Let ${\cal M}$
be a maximal set in $\mathbb{D}^{m\times n}$. Then $\ell$ is a line in $AG({\cal M})$ if and only if
$\ell=\mathcal{M}\cap \mathcal{N}$, where $\mathcal{M}$ and $\mathcal{N}$ are two  maximal sets of different types with $\mathcal{M}\cap\mathcal{N}\neq\emptyset$.
Moreover,  $|\ell|=|\mathcal{M}\cap\mathcal{N}|=|\mathbb{D}|$.
\end{lemma}

\begin{lemma}\label{non-degeneratelemma00a}  Let $m,n,m', n'\geq 2$ be integers. Suppose that
$\varphi:  {\mathbb{D}}^{m\times n}\rightarrow  {\mathbb{D}'}^{m'\times n'}$ is a non-degenerate graph homomorphism with $\varphi(0)=0$. Then
\begin{itemize}
\item[{\em (a)}]   if ${\cal M}$ is a maximal set of type one [resp. type two]  containing $0$ in ${\mathbb{D}}^{m\times n}$ and $\varphi({\cal M})\subseteq {\cal M}'$, where
 ${\cal M}'$ is a maximal set  in ${\mathbb{D}'}^{m'\times n'}$, then for any $A\in{\cal M}$,  there exists a  maximal set ${\cal R}$
 of type one [resp. type two]  in ${\mathbb{D}}^{m\times n}$, such that $\varphi({\cal R})\subseteq {\cal R}'$ and ${\cal R}\cap {\cal M}=\{A\}$,
where ${\cal R}'$ is a maximal set  in ${\mathbb{D}'}^{m'\times n'}$, ${\cal R}'$ and ${\cal M}'$ are of the same type with ${\cal R}'\neq {\cal M}'$;

\item[{\em (b)}] if ${\cal M}$ and ${\cal N}$ are two distinct maximal sets of the same type  [resp. different types]
in ${\mathbb{D}}^{m\times n}$ such that $0\in \mathcal{M}$ and ${\cal M}\cap {\cal N}\neq \emptyset$, then  $\varphi({\cal M})\subseteq {\cal M}'$
and $\varphi({\cal N})\subseteq {\cal N}'$,
where ${\cal M}'$ and ${\cal N}'$ are two maximal sets of the same type  [resp. different types]
in ${\mathbb{D}'}^{m'\times n'}$;

\item[{\em (c)}] if ${\cal M}$ is a maximal set containing $0$ in ${\mathbb{D}}^{m\times n}$  and  $\varphi(\mathcal{M})\subseteq \mathcal{M}'$
where ${\cal M}'$ is a maximal set in ${\mathbb{D}'}^{m'\times n'}$,
then ${\cal M}'$ is the unique  maximal set containing $\varphi({\cal M})$, and $\varphi({\cal M})$ is not contained in any line of $AG({\cal M}')$.
\end{itemize}
\end{lemma}
\proof
(a).    Let ${\cal M}$ be a maximal set of type one [resp. type two] containing $0$ in ${\mathbb{D}}^{m\times n}$, and let $\varphi({\cal M})\subseteq {\cal M}'$ where
 ${\cal M}'$ is a maximal set  in ${\mathbb{D}'}^{m'\times n'}$.
 Without loss of generality we assume that both $\mathcal{M}$ and $\mathcal{M}'$ are of type one.

Let  matrix $A\in{\cal M}$. Then $A\in \mathbb{D}_{\leq 1}^{m\times n}$.
By Lemma \ref{Rectangular-PID2-4}, ${\cal M}=P{\cal M}_1+A$ where $P\in GL_m(\mathbb{D})$.
There exists a type one maximal set ${\cal R}$  in ${\mathbb{D}}^{m\times n}$,
such that ${\cal R}\neq {\cal M}$, ${\cal R}\cap {\cal M}=\{A\}$ and $\varphi({\cal R})\subseteq {\cal R}'$, where ${\cal R}'$ is a  maximal set in ${\mathbb{D}'}^{m'\times n'}$
and ${\cal R}'\neq {\cal M}'$.
 Otherwise, for any type one maximal set
 ${\cal R}$ in ${\mathbb{D}}^{m\times n}$ containing $A$, we have $\varphi({\cal R})\subseteq {\cal M}'$.
Since every matrix in $\mathbb{D}_{\leq 1}^{m\times n}+A$ is contained in some type one maximal set which contais $A$,
we get $\varphi\left(\mathbb{D}_{\leq 1}^{m\times n}+A\right)\subseteq {\cal M}'$,
a contradiction to $\varphi$ being non-degenerate.

We show that  ${\cal R}'$ is of type one by contradiction. Suppose that ${\cal R}'$ is of type two.
Let ${\cal N}$ be a maximal set of type two in ${\mathbb{D}}^{m\times n}$ containing $A$, and let $\varphi({\cal N})\subseteq {\cal N}'$
where ${\cal N}'$ is a  maximal set in ${\mathbb{D}'}^{m'\times n'}$. By Corollary \ref{Rectangular-1-13}, we have
 $|{\cal N}\cap {\cal R}|\geq 2$ and $|{\cal N}\cap {\cal M}|\geq 2$,    hence
$|{\cal N}'\cap {\cal R}'|\geq 2$ and $|{\cal N}'\cap {\cal M}'|\geq 2$. Applying Corollary \ref{Rectangular-1-13} again, we get either
${\cal N}'={\cal R}'$ or ${\cal N}'={\cal M}'$, and hence ${\cal N}'\subseteq {\cal M}'\cup {\cal R}'$.
It follows that  $\varphi({\cal N})\subseteq {\cal M}'\cup {\cal R}'$ for any type two maximal set ${\cal N}$ in ${\mathbb{D}}^{m\times n}$ containing $A$.
Since every matrix in $\mathbb{D}_{\leq 1}^{m\times n}+A$ is contained in some type two maximal set containing $A$,
we get $\varphi\left(\mathbb{D}_{\leq 1}^{m\times n}+A\right)\subseteq {\cal M}'\cup {\cal R}'$,
a contradiction to $\varphi$ being non-degenerate. Thus,  ${\cal R}'$ is of type one.

(b). \ Let ${\cal M}$ and ${\cal N}$ be two distinct maximal sets of the same type  [resp. different types]
in ${\mathbb{D}}^{m\times n}$ such that $0\in \mathcal{M}$ and ${\cal M}\cap {\cal N}\neq \emptyset$.
 Then there are  maximal sets ${\cal M}'$ and ${\cal N}'$ in ${\mathbb{D}'}^{m'\times n'}$,
such that $\varphi({\cal M})\subseteq {\cal M}'$ and $\varphi({\cal N})\subseteq {\cal N}'$.
Without loss of generality, we assume that ${\cal M}$ is of type one. Let $A\in{\cal M}\cap {\cal N}$.
Then $A\in \mathbb{D}_{\leq 1}^{m\times n}$ because $0\in\mathcal{M}$.
We prove that  ${\cal M}'$ and ${\cal N}'$ are of the same type  [resp. different types] as follows.

{\em Case 1.} \ ${\cal M}$ and ${\cal N}$ are  of the same type (i.e., type one).
By the (a) and Lemma \ref{Rectangular-PID2-4}, there exist  two distinct type two maximal sets ${\cal L}$ and ${\cal T}$ in
${\mathbb{D}}^{m\times n}$, such that  $0\in\mathcal{L}$, $\mathcal{L}\cap \mathcal{T}=\{A\}$, $\varphi({\cal L})\subseteq {\cal L}'$ and  $\varphi({\cal T})\subseteq {\cal T}'$,
where ${\cal L}'$ and ${\cal T}'$ are two distinct maximal sets of the same type.

Since $|{\cal N}\cap {\cal L}|\geq 2$ and $|{\cal N}\cap {\cal T}|\geq 2$, we have
$|{\cal N}'\cap {\cal L}'|\geq 2$ and $|{\cal N}'\cap {\cal T}'|\geq 2$. By Corollary \ref{Rectangular-1-13} and ${\cal L}'\neq {\cal T}'$, we get that
 ${\cal N}'$ and ${\cal L}'$ are of  different types. Similarly,  ${\cal M}'$ and ${\cal L}'$ are of  different types.
 Consequently,  ${\cal M}'$ and ${\cal N}'$ are of the same type.

{\em Case 2.} \ ${\cal M}$ and ${\cal N}$ are  of different types. Then ${\cal N}$ is of type two. By the (a),  there exists a type one maximal set
${\cal R}$ in ${\mathbb{D}}^{m\times n}$ containing $A$, such that $\varphi({\cal R})\subseteq {\cal R}'$ and ${\cal R}\neq {\cal M}$,
where ${\cal R}'$ is a maximal set in ${\mathbb{D}'}^{m'\times n'}$ such that ${\cal R}'$ and ${\cal M}'$ are  of the same type and ${\cal R}'\neq {\cal M}'$.
Since $|{\cal N}\cap {\cal R}|\geq 2$ and $|{\cal N}\cap {\cal M}|\geq 2$, one has
$|{\cal N}'\cap {\cal R}'|\geq 2$ and $|{\cal N}'\cap {\cal M}'|\geq 2$. By Corollary \ref{Rectangular-1-13}, we obtain that
 ${\cal N}'$ and ${\cal M}'$ are of  different types.

Combining Case 1 with Case 2,  the (b) of this lemma holds.

(c). \ Let ${\cal M}$ be a maximal set containing $0$ in ${\mathbb{D}}^{m\times n}$,  and  let $\varphi(\mathcal{M})\subseteq \mathcal{M}'$
where ${\cal M}'$ is a maximal set in ${\mathbb{D}'}^{m'\times n'}$. Without loss of generality,
we can assume that both ${\cal M}$ and ${\cal M}'$ are of type one.
 By the (a) and  Lemma \ref{Rectangular-PID2-4}, there exist  two distinct type two maximal sets ${\cal L}$ and ${\cal T}$ in
${\mathbb{D}}^{m\times n}$, such that ${\cal M}\cap{\cal L}\cap{\cal T}=\{0\}$,   $\varphi({\cal L})\subseteq {\cal L}'$ and  $\varphi({\cal T})\subseteq {\cal T}'$,
where ${\cal L}'$ and ${\cal T}'$ are  two distinct maximal sets of the same type.

Since $|{\cal M}\cap {\cal L}|\geq 2$ and $|{\cal M}\cap {\cal T}|\geq 2$, we have $|{\cal M}'\cap {\cal L}'|\geq 2$ and $|{\cal M}'\cap {\cal T}'|\geq 2$.
Thus  Corollary \ref{Rectangular-1-13} implies that
both ${\cal L}'$ and ${\cal T}'$ are of type two.
By Lemma \ref{Rectangular-PID2-11}, $\ell_1:=\mathcal{M}\cap \mathcal{L}$ and $\ell_2:=\mathcal{M}\cap \mathcal{T}$ are two lines in $AG({\cal M})$;
 $\ell_1':=\mathcal{M}'\cap \mathcal{L}'$ and $\ell_2':=\mathcal{M}'\cap \mathcal{T}'$ are also two lines in $AG({\cal M}')$.
Clearly, $\varphi(\ell_i)\subseteq \ell_i'$, $i=1,2$, and $\ell'_1\neq \ell'_2$.
Since two different lines have at most a common point,
it follows that $\varphi({\cal M})$ is not contained in any line in $AG({\cal M}')$.
Applying Lemma \ref{Rectangular-PID2-11} and Corollary \ref{Rectangular-1-13}, it is clear that ${\cal M}'$ is the unique  maximal set containing $\varphi({\cal M})$.
$\qed$

If a map $\varphi:  {\mathbb{D}}^{m\times n}\rightarrow  {\mathbb{D}'}^{m'\times n'}$ preserves distances $1$ and $2$, then $\varphi$
carries  distinct maximal sets in ${\mathbb{D}}^{m\times n}$  into  distinct maximal sets in ${\mathbb{D}'}^{m'\times n'}$,
and hence $\varphi$ is a non-degenerate graph homomorphism.
Thus,  by Lemma \ref{non-degeneratelemma00a}(b) we have the following corollary.

\begin{corollary}\label{degenerate-aa4}
Let $\varphi:  {\mathbb{D}}^{m\times n}\rightarrow  {\mathbb{D}'}^{m'\times n'}$ be a distances $1$ and $2$ preserving map with $\varphi(0)=0$.
Suppose that  ${\cal M}$ and ${\cal N}$ are two distinct maximal sets of the same type  [resp. different types]
in ${\mathbb{D}}^{m\times n}$ with $0\in\mathcal{M}$ and ${\cal M}\cap {\cal N}\neq \emptyset$. Then there are two distinct maximal sets ${\cal M}'$
and ${\cal N}'$ in ${\mathbb{D}'}^{m'\times n'}$, such that $\varphi({\cal M})\subseteq {\cal M}'$, $\varphi({\cal N})\subseteq {\cal N}'$,
and  ${\cal M}', {\cal N}'$ are of the same type  [resp. different types].
\end{corollary}

Let $V$ be a left   vector space over  $\mathbb{D}$, and let $V'$ be a  left [resp. right]  vector space over  $\mathbb{D}'$.
Suppose $f: V\rightarrow V'$ is a map. The map $f$ is called a {\em semi-affine map} if there exists a ring homomorphism  $\tau$ [resp. ring anti-homomorphism $\sigma$]
from $\mathbb{D}$ to $\mathbb{D}'$, such that
$f(\alpha x+\beta y+\gamma z)=\alpha^\tau f(x)+\beta^\tau f(y)+\gamma^\tau f(z)$
[resp. $f(\alpha x+\beta y+\gamma z)= f(x)\alpha^\sigma+ f(y)\beta^\sigma+ f(z)\gamma^\sigma]$,
for all $x, y, z\in V$ and $\alpha, \beta, \gamma\in \mathbb{D}$ with $\alpha+\beta+\gamma=1$. A semi-affine map $f$ is called a {\em quasi-affine map}
 if  $\tau$ [resp. $\sigma$] is a ring isomorphism [resp. ring anti-isomorphism].  A quasi-affine map $f$ is called an {\em affine map}
if $\mathbb{D}=\mathbb{D}'$ and $\tau$ [resp. $\sigma$] is the identity (cf. \cite[Definition 12.5.1]{ModernProjective}).

For two right  vector spaces  over division rings, or a right  vector space and a left  vector space  over division rings,
we can define similarly the semi-affine map,  quasi-affine map and affine map between two vector spaces, respectively.

Suppose  both $V$ and $V'$ are left [resp. right] vector spaces over $\mathbb{D}$ and $\mathbb{D}'$, respectively. An additive map $f:
V\rightarrow V'$ is called a {\em semi-linear map} if  there exists a ring homomorphism $\tau$ from $\mathbb{D}$ to $\mathbb{D}'$ such that
$f(\lambda x)=\lambda^\tau f(x)$  [resp. $f(x\lambda)= f(x)\lambda^\tau$],  for all $x\in V$, $\lambda\in \mathbb{D}$.
For  a left [resp. right] vector space $V$ over $\mathbb{D}$ and a  right  [resp. left] vector space $V'$ over $\mathbb{D}'$,
we can define similarly a semi-linear map from $V$ to $V'$.
A semi-linear map is called {\em quasi-linear map} if $\tau$ [resp. $\sigma$] is
 a ring isomorphism [resp. ring anti-isomorphism]. A quasi-linear map is called {\em linear map} if $V=V'$ and $\tau$ [resp. $\sigma$] is the identity.

\begin{remark}\label{weaksemilinearandquasilinear2}
Let $\mathbb{D},\mathbb{D}'$ be division rings, and let $n, m$ be positive integers with $n\geq 2$. Then
\begin{itemize}
\item[{\rm (a)}]  if  $V={\mathbb{D}}^n$ $[$resp. $V= \,^n\mathbb{D}$$]$,  $V'=\mathbb{D}^{m}$ $[$resp.  $V'= \,^{m}\mathbb{D}$$]$  and
$f: V\rightarrow V'$ is a  semi-linear map, then $f$ is of the form
$f(x)=x^\tau P$ [resp. $f(x)=Q x^\tau$],  $x\in V$,
where $\tau$  is a ring homomorphism  from $\mathbb{D}$ to $\mathbb{D}'$, and $P\in {\mathbb{D}'}^{n\times m}$ $[$resp.  $Q\in {\mathbb{D}'}^{m\times n}$$]$ is fixed;

\item[{\rm (b)}]  if  $V={\mathbb{D}}^n$ $[$resp. $V= \,^n\mathbb{D}$$]$,  $V'= \, ^m\mathbb{D}$ $[$resp.  $V'=\mathbb{D}^{m}$$]$ and
$f: V\rightarrow V'$ is a  semi-linear map, then $f$ is of the form
$f(x)=P\, ^tx^\sigma$ [resp. $f(x)=\,^tx^\sigma Q$],  $x\in V$,
where  $\sigma$ is a ring anti-homomorphism from $\mathbb{D}$ to $\mathbb{D}'$, and $P\in {\mathbb{D}'}^{m\times n}$ $[$resp.  $Q\in {\mathbb{D}'}^{n\times m}$$]$ is  fixed.
\end{itemize}
\end{remark}

\begin{lemma}{\em (see \cite[Proposition 12.5.2]{ModernProjective})}\label{semi-affine map-a-001} \
Let $V$ and $V'$ be left or right vector spaces  over  division rings $\mathbb{D}$ and $\mathbb{D}'$, respectively. Assume that $f: V\rightarrow V'$ is a map. Then
\begin{itemize}
\item[\em(a)] $f$ is a semi-affine map if and only if  $g(x):=f(x)-f(0)$ is a  semi-linear map;

\item[\em(b)] $f$ is a quasi-affine map if and only if  $g(x):=f(x)-f(0)$ is a  quasi-linear map;

\item[\em(c)] $f$ is an affine map if and only if  $g(x):=f(x)-f(0)$ is a  linear map.
\end{itemize}
\end{lemma}

Let $V$ and $V'$ be  left or right  vector spaces of dimensions $\geq 2$ over division rings. A  map $f: V\rightarrow V'$ is called an {\em affine morphism} if
$a\in b\vee c$ $\Rightarrow$ $f(a)\in f(b)\vee f(c)$ for all $a,b,c\in V$ (i.e., $f$ maps lines into lines).
An affine morphism is called {\em non-degenerate} if its image is not contained in a line.

The following definition is a simpler version of the weighted semi-affine map \cite[Definitian 12.5.6]{ModernProjective},
since we do not need to use ``partial map" of \cite{ModernProjective}.

Suppose  $V$ is a left  vector space over a division ring $\mathbb{D}$, $V'$ a left [resp. right]  vector space over a division ring $\mathbb{D}'$,
and their dimensions $\geq 2$.
A  map $g: V\rightarrow V'$ is called a {\em weighted semi-affine map} if there exists a semi-affine map
$(k, h):   V\rightarrow \mathbb{D}'\times V'$ [resp.  $(h, k):   V\rightarrow V'\times \mathbb{D}'$], where
$h$ is a  semi-affine map from $V$ to $V'$ and $k$ is a  semi-affine map from $V$ to $\mathbb{D}'$,  such that $k(x)\neq 0$ for all $x\in V$ and
\begin{equation}\label{weightedsemi-affine1.1}
g(x)=(k(x))^{-1}h(x) \ \ [{\rm resp.} \ g(x)=h(x)(k(x))^{-1}],   \  x\in V.
\end{equation}
Note that in  $(k, h)$ [resp. $(h,k)$], semi-affine maps $h$ and $k$ have the same ring homomorphism [resp. ring anti-homomorphism] from $\mathbb{D}$ to $\mathbb{D}'$
associated to $h$ and $k$.

For two right  vector spaces  over division rings, or a right  vector space  and a left  vector space over division rings,
we can define similarly the weighted semi-affine map between two vector spaces.
Every weighted semi-affine map  is an affine morphism (see \cite[Proposition 12.5.7]{ModernProjective}).

\begin{lemma}\label{weaksemilinea34ear2}
Let $\mathbb{D}, \mathbb{D}'$ be division rings, and let $m, n\geq 2$ be integers. Then:
\begin{itemize}
\item[{\em (a)}]
if  $V={\mathbb{D}}^n$,  $V'= {\mathbb{D}'}^m$ $[$resp. $V= \,^n\mathbb{D}$, $V'= \,^m\mathbb{D}'$$]$ and
 $g: V\rightarrow V'$ is a  weighted semi-affine map with $g(0)=0$, then $g$ is of the form
$$
g(x)=(k(x))^{-1}x^\tau P \ \ [resp. \ g(x)=Qx^\tau (k'(x))^{-1}],  \  x\in V,
$$
where $\tau$ is a ring homomorphism  from $\mathbb{D}$ to $\mathbb{D}'$,  $P\in {\mathbb{D}'}^{n\times m}$ $[$resp.  $Q\in {\mathbb{D}'}^{m\times n}$$]$
is fixed, and $k(x)$ is of the form $k(x)=\sum_{i=1}^nx_i^\tau a_i+b$ with $k(x)\neq 0$ for all $x=(x_1,\ldots, x_n)\in V$
$[$resp. $k'(x)=\sum_{i=1}^na_ix_i^\tau +b$ with $k'(x)\neq 0$ for all $x=\,^t(x_1,\ldots, x_n)\in V$$]$;

\item[{\em (b)}]  if  $V={\mathbb{D}}^n$,  $V'= \,^m\mathbb{D}'$ $[$resp.  $V= \,^n\mathbb{D}$,  $V'={\mathbb{D}'}^m$$]$ and
 $g: V\rightarrow V'$ is a  weighted semi-affine map with $g(0)=0$, then $g$ is of the form
$$
g(x)=P\,{^tx^\sigma} (k(x))^{-1} \ \ [resp. \ g(x)= (k'(x))^{-1}\,^tx^\sigma Q],  \  x\in V,
$$
where $\sigma$ is a ring anti-homomorphism from $\mathbb{D}$ to $\mathbb{D}'$, $P\in {\mathbb{D}'}^{m\times n}$ $[$resp.  $Q\in {\mathbb{D}'}^{n\times m}$$]$ is fixed, and
$k(x)$ is of the form  $k(x)=\sum_{i=1}^na_ix_i^\sigma +b$ with $k(x)\neq 0$ for all $x=(x_1,\ldots, x_n)\in V$ $[$resp. $k'(x)=\sum_{i=1}^nx_i^\sigma a_i +b$
with $k'(x)\neq 0$ for all $x= \,^t(x_1,\ldots, x_n)\in V$$]$.
\end{itemize}
\end{lemma}
\proof(a). Suppose that  $V={\mathbb{D}}^n$  and $V'= {\mathbb{D}'}^m$.
 Let $g: V\rightarrow V'$ be a  weighted semi-affine map with $g(0)=0$.
By the definition of the weighted semi-affine map, there exists a semi-affine map $(k, h):   V\rightarrow \mathbb{D}'\times V'$, where
$h$ is a  semi-affine map from $V$ to $V'$ and $k$ is a  semi-affine map from $V$ to $\mathbb{D}'$,  such that $k(x)\neq 0$  and
$$
g(x)=(k(x))^{-1}h(x),   \ x\in V.
$$
Let $\tau: \mathbb{D}\rightarrow \mathbb{D}'$ be the ring homomorphism  associated to the semi-affine map $(k,h)$.

By Lemma \ref{semi-affine map-a-001}, $k'(x):=k(x)-k(0)$ is a semi-linear map from $\mathbb{D}^n$ to $\mathbb{D}'$,
and  $h'(x):=h(x)-h(0)$ is a semi-linear map from $\mathbb{D}^n$ to ${\mathbb{D}'}^m$. Since $g(0)=0$, we get $h(0)=0$. Hence
$h(x)=h'(x)$ is a semi-linear map.
Let $e_i$ is the $i$-th row of $I_n$, $i=1, \ldots, n$. For any $x=\sum_{i=1}^nx_ie_i\in V$ where $x_i\in \mathbb{D}$,
we have $k'(x)=k'(\sum_{i=1}^nx_ie_i)=\sum_{i=1}^nx_i^\tau k'(e_i)$.
Thus $k(x)=\sum_{i=1}^nx_i^\tau k'(e_i)+k(0)$. Let $a_i=k'(e_i)$ and $b=k(0)$. Then  $k(x)=\sum_{i=1}^nx_i^\tau a_i+b$.
By Remark \ref{weaksemilinearandquasilinear2}, $h(x)=x^\tau P$ where $P\in {\mathbb{D}'}^{n\times m}$ is fixed. Therefore,
$$g(x)=(k(x))^{-1}x^\tau P,  \   x\in V.
$$

For the case of  $V= \,^n\mathbb{D}$ and  $V'= \,^m\mathbb{D}'$, we have similarly that $g(x)=Qx^\tau (k'(x))^{-1}$ ($x\in V$),
where $Q\in {\mathbb{D}'}^{m\times n}$ is fixed, and $k'(x)=\sum_{i=1}^na_ix_i^\tau +b$ with $k'(x)\neq 0$ for all $x=\,^t(x_1,\ldots, x_n)\in V$.

(b). The proof is similar.
$\qed$

\begin{proposition}{\em(see \cite[Proposition 12.5.9]{ModernProjective})}\label{weighted semi-affine3} \
 Suppose  $V$ and $V$ are left vector spaces of  dimensions $\geq 2$ over  division rings $\mathbb{D}$ and $\mathbb{D}'$, respectively.
Let $g: V\rightarrow V'$ is a non-constant weighted semi-affine map  which is defined by the semi-affine map $(k, h): V\rightarrow  \mathbb{D}'\times V'$. Then the
following three conditions are equivalent:
\begin{itemize}
\item[\em(a)] $g$ is a quasi-affine map (with empty kernel),

\item[\em(b)] there exist $a,b\in V$ such that $g(a)\neq g(b)$ and $g(a\vee b)=g(a)\vee g(b)$,

\item[\em(c)] the ring homomorphism $\tau: \mathbb{D}\rightarrow \mathbb{D}'$ associated to $(k, h)$ is an isomorphism.
\end{itemize}
\end{proposition}

We have the following fundamental theorem of the affine geometry on vector spaces.

\begin{theorem}\label{weighted semi-affine9}{\em (see \cite[Theorem 12.6.7]{ModernProjective})} \
Let $\mathbb{D}, \mathbb{D}'$ be division rings with $|\mathbb{D}|\geq 4$. Suppose that
 $V$ and  $V'$ are left or right vector spaces of dimensions $\geq 2$ over  $\mathbb{D}$ and $\mathbb{D}'$, respectively. Then every
non-degenerate affine morphism  $f: V\rightarrow V'$ is a weighted semi-affine map.
\end{theorem}

Let $\varphi:  {\mathbb{D}}^{m\times n}\rightarrow  {\mathbb{D}'}^{m'\times n'}$ be a graph homomorphism with $\varphi(0)=0$.
Suppose  $\mathcal{M}$ is a  maximal set of type one containing $0$ in ${\mathbb{D}}^{m\times n}$ and
$\varphi({\cal M})\subseteq {\cal M}'$, where  $\mathcal{M}'$ is a  maximal set of type one containing $0$ in ${\mathbb{D}'}^{m'\times n'}$.
Then there are $P\in GL_m(\mathbb{D})$ and $P'\in GL_{m'}(\mathbb{D})$, such that
$\mathcal{M}=P\mathcal{M}_1$ and $\mathcal{M}'=P'\mathcal{M}_1'$. Let
$$\varphi\left(P\left(
         \begin{array}{c}
           x\\
           0 \\
         \end{array}
       \right)\right)=
P'\left(\begin{array}{c}
x^\rho \\
0 \\
\end{array}
\right), \ x\in {\mathbb{D}}^n,$$
 where $\rho: \mathbb{D}\rightarrow \mathbb{D}'$ is an injective map with $0^\rho=0$.
Then the restriction map $\varphi\mid_{\mathcal{M}}: \mathcal{M}\rightarrow \mathcal{M}'$ induces the map $\varphi': {\mathbb{D}}^n\rightarrow {\mathbb{D}'}^{n'}$
by $\varphi'(x)=x^\rho$ for all $x\in {\mathbb{D}}^n$.
 The restriction map $\varphi\mid_{\mathcal{M}}$
 is called  a {\em weighted semi-affine map} [resp. {\em non-degenerate affine morphism}] if the map  $\varphi'$ is a weighted semi-affine map [resp.  non-degenerate affine morphism].

 Similarly, when  $\varphi({\cal M})\subseteq {\cal N}'$ where ${\cal N}'$ is a maximal set of type two containing $0$, or
 ${\cal M}$ is a  maximal set of type two containing $0$, we can define similarly the weighted semi-affine map (or non-degenerate affine morphism)
 for the restriction map $\varphi\mid_{\mathcal{M}}$.

\begin{lemma}\label{degenerate-2} Let $\mathbb{D}, \mathbb{D}'$ be division rings with $|\mathbb{D}|\geq 4$, and
let $m,n,m', n'\geq 2$ be integers. Suppose that $\varphi:  {\mathbb{D}}^{m\times n}\rightarrow  {\mathbb{D}'}^{m'\times n'}$
is a non-degenerate graph homomorphism with $\varphi(0)=0$. Let ${\cal M}$ be a maximal set containing $0$ in ${\mathbb{D}}^{m\times n}$, such that
$\varphi({\cal M})\subseteq {\cal M}'$ where ${\cal M}'$ is a maximal set containing $0$ in ${\mathbb{D}'}^{m'\times n'}$.
Then the restriction map $\varphi\mid_{\mathcal{M}}: \mathcal{M}\rightarrow \mathcal{M}'$
 is an injective weighted semi-affine map. Moreover, if $\mathcal{S}$ is an adjacent set in ${\mathbb{D}}^{m\times n}$ with  $0\in \mathcal{S}$
 and $|\mathcal{S}|\geq 2$, then $${\rm dim}(\varphi(\mathcal{S}))\leq {\rm dim}(\mathcal{S}).$$
\end{lemma}
\proof
Without loss of generality, we assume that both  ${\cal M}$ and  ${\cal M}'$ are   maximal sets of type one.
Let $\ell$ be a line in $AG({\cal M})$. By Lemma \ref{Rectangular-PID2-11}, $\ell={\cal M}\cap{\cal N}$, where ${\cal N}$ is a maximal set of type two
 in ${\mathbb{D}}^{m\times n}$ with ${\cal M}\cap{\cal N}\neq \emptyset$. There exists a  maximal set ${\cal N}'$  in ${\mathbb{D}'}^{m'\times n'}$ such that
 $\varphi(\mathcal{N})\subseteq {\cal N}'$ and $|{\cal M}'\cap {\cal N}'|\geq 2$. By Lemma \ref{non-degeneratelemma00a}(b),
 ${\cal N}'$ is of type two. Hence Lemma \ref{Rectangular-PID2-11} implies that  $\ell'={\cal M}'\cap {\cal N}'$ is a
 line in $AG(\mathcal{M}')$ and $\varphi(\ell)\subseteq \ell'$. Thus, $\varphi$ maps  lines in $AG(\mathcal{M})$ into  lines in $AG(\mathcal{M}')$.
 By Lemma \ref{non-degeneratelemma00a}(c), $\varphi({\cal M})$ is not contained in any line in $AG({\cal M}')$.
Therefore,  the restriction map $\varphi\mid_{\mathcal{M}}: \mathcal{M}\rightarrow \mathcal{M}'$ is a non-degenerate affine morphism.
 By Theorem \ref{weighted semi-affine9}, the restriction map $\varphi\mid_{\mathcal{M}}$ is a weighted semi-affine map.
 Clearly, $\varphi\mid_{\mathcal{M}}$  is injective.

Let $\mathcal{S}$ be an adjacent set in ${\mathbb{D}}^{m\times n}$ with $0\in \mathcal{S}$ and $|\mathcal{S}|\geq 2$.
By Lemma \ref{Rectangular-PID2-4},we have either $\mathcal{S}\subseteq P\mathcal{M}_1$
or $\mathcal{S}\subseteq\mathcal{N}_1Q$, where $P\in GL_m(\mathbb{D})$, $Q\in GL_n(\mathbb{D})$.
Thus, $\varphi(\mathcal{S})\subseteq \varphi(P\mathcal{M}_1)$ or $\varphi(\mathcal{S})\subseteq \varphi(\mathcal{N}_1Q)$.
Without loss of generality, we assume that $\mathcal{S}\subseteq P\mathcal{M}_1$ and $\varphi(P\mathcal{M}_1)\subseteq P'\mathcal{M}_1'$,
where $P'\in GL_{m'}(\mathbb{D}')$. Then $\varphi(\mathcal{S})\subseteq P'\mathcal{M}_1'$ and
the restriction map $\varphi\mid_{P\mathcal{M}_1}: P\mathcal{M}_1\rightarrow P'\mathcal{M}_1'$
is a weighted semi-affine map. By Lemma \ref{weaksemilinea34ear2}, $\varphi\mid_{P\mathcal{M}_1}$ is of the form
$$
\varphi\left(P\left(
         \begin{array}{c}
           x \\
           0 \\
         \end{array}
       \right)\right)=
P'\left(\begin{array}{c}
(k(x))^{-1}x^{\tau}T \\
0 \\
\end{array}
\right),  \  x\in {\mathbb{D}}^n,
$$
where  $\tau: \mathbb{D}\rightarrow \mathbb{D}'$ is a nonzero ring homomorphism, $T\in {\mathbb{D}'}^{n\times n'}$, and $k(x)$ is of the form
$k(x)=\sum_{j=1}^nx_{1j}^{\tau} a_{j1}+b\neq 0$ for all $x=(x_{11}, \ldots, x_{1n})\in {\mathbb{D}}^n$.
Thus, it is easy to see that ${\rm dim}(\varphi(\mathcal{S}))\leq {\rm dim}(\mathcal{S})$.
$\qed$

\section{Grassmann Spaces and Related Lemmas}

In this section, we  introduce  Grassmann spaces (cf. \cite{W3,Huang-book, L.P.Huang2009-3}).

 Let  $n\geq 3$ be an integer, and let $\mathscr{P}(\mathbb{D}^{n})$ [resp. $\mathscr{P}(^{n}\mathbb{D})$] be the $(n-1)$-dimensional
 left [resp. right] {\em projective geometry} on $\mathbb{D}^{n}$ [resp. $^{n}\mathbb{D}$]. In  $\mathscr{P}(\mathbb{D}^{n})$ [reap. $\mathscr{P}(^{n}\mathbb{D})$],
 an $(r+1)$-dimensional vector subspace $W:=[v_1, \ldots, v_{r+1}]$
 is called an $r$-dimensional {\em projective flat} ({\em $r$-flat} for short), and the matrix
$$\mbox{$\left(
\begin{array}{c}
 v_1 \\
\vdots\\
v_{r+1} \\
\end{array}\right)$ \ [resp. $(v_1, \ldots, v_{r+1})$]}$$
is called a {\em matrix representation} of the $r$-flat $W$. For simpleness, the matrix representation of an $r$-flat $W$ is also denoted
by $W$.

Let  $n,m\geq 2$ be integers. The $(m-1)$-dimensional left [resp. right] {\em Grassmann space} over $\mathbb{D}$,  denoted
 by $\mathscr{G}^l_{m+n-1,m-1}(\mathbb{D})$ [resp.  $\mathscr{G}^r_{m+n-1,m-1}(\mathbb{D})$], is the set of all $(m-1)$-flats
 in $\mathscr{P}(\mathbb{D}^{m+n})$ [resp. $\mathscr{P}(^{m+n}\mathbb{D})$].
 In a Grassmann space, an $(m-1)$-flat $W$ is called a {\em point}. Matrices $W$ and $W'$ of rank $m$ are two matrix representations
of one point in  $\mathscr{G}^l_{m+n-1,m-1}(\mathbb{D})$ [resp.  $\mathscr{G}^r_{m+n-1,m-1}(\mathbb{D})$]
if and only if there is a $G\in GL_m(\mathbb{D})$ such that $W'=GW$ [resp. $W'=WG$].

 In the Grassmann space $\mathscr{G}^l_{m+n-1,m-1}(\mathbb{D})$ [resp.  $\mathscr{G}^r_{m+n-1,m-1}(\mathbb{D})$], two points $W_1$ and $W_2$
are said to be of {\em arithmetic distance} $r$, denoted by {\rm ad}$(W_1,W_2)=r$, if  $W_1\cap W_2$  is an $(m-r-1)$-flat.
 Two points $W_1$ and $W_2$ are called {\em adjacent} if {\rm ad}$(W_1,W_2)=1$. For matrix representations $W_1$ and $W_2$ of points $W_1$ and $W_2$,
 we have (cf. \cite{W3,Huang-book})
\begin{equation}\label{Grassmanndistance}
\mbox{${\rm ad}(W_1,W_2)={\rm rank}\left(\begin{array}{c}
  W_1 \\
  W_2 \\
\end{array}\right)-m$ \ [resp. ${\rm ad}(W_1,W_2)={\rm rank}(W_1, W_2)-m$].}
\end{equation}

Recall that the projective general linear group of degree $m+n$ over $\mathbb{D}$ is $PGL_{m+n}(\mathbb{D})=GL_{m+n}(\mathbb{D})/Z_{m+n}$, where
$Z$ is the center of $\mathbb{D}$ and $Z_{m+n}=\{\lambda I_{m+n}: \lambda \in Z^*\}$. The $PGL_{m+n}(\mathbb{D})$ can be regarded as a group of motions
on $\mathscr{G}^l_{m+n-1,m-1}(\mathbb{D})$  [resp. $\mathscr{G}^r_{m+n-1,m-1}(\mathbb{D})$]. That is, each element $T\in PGL_{m+n}(\mathbb{D})$ defines a bijective map
$W\mapsto  WT$ [resp. $W\mapsto  TW$] from $\mathscr{G}^l_{m+n-1,m-1}(\mathbb{D})$ [resp. $\mathscr{G}^r_{m+n-1,m-1}(\mathbb{D})$] to itself.
By   \cite[Propositions 3.28 and 3.31]{W3}, the group $PGL_{m+n}(\mathbb{D})$ acts transitively on  $\mathscr{G}^l_{m+n-1,m-1}(\mathbb{D})$ [resp. $\mathscr{G}^r_{m+n-1,m-1}(\mathbb{D})$],
and every element of the group $PGL_{m+n}(\mathbb{D})$  preserves  the arithmetic distance between any two points.

Suppose that  $m, n\geq 2$ are integers and $k, r$ are integers with $1\leq r\leq k\leq {\rm max}\{m,n\}$.
 We define the following mathematical symbols.
\begin{itemize}
\item
Let $\mathcal{M}_{k,r}(\mathbb{D}^{m\times n})$ denote the set of all matrices $X\in \mathbb{D}_r^{m\times n}$
having exactly $k$ nonzero rows.  Write $\alpha=\{i_1,i_2, \ldots, i_k\}$ where $1\leq i_1<i_2< \cdots <i_k\leq m$.
Let $\mathcal{M}_{k,r}^\alpha(\mathbb{D}^{m\times n})$ denote the set of all matrices in $\mathbb{D}_r^{m\times n}$ that
have exactly $k$ nonzero rows  being row $i_1$, row $i_2$,  $\ldots$, row $i_k$.

\item  Let $\mathcal{N}_{k,r}(\mathbb{D}^{m\times n})$  denote the set of all matrices $X\in \mathbb{D}_r^{m\times n}$
having exactly $k$ nonzero columns. Write $\beta=\{j_1,j_2, \ldots, j_k\}$ where $1\leq j_1<j_2< \cdots <j_k\leq n$.
 Let $\mathcal{N}_{k,r}^\beta(\mathbb{D}^{m\times n})$  denote the set of all matrices in $\mathbb{D}_r^{m\times n}$ that
have exactly $k$ nonzero columns  being column $j_1$, column $j_2$,  $\ldots$, column $j_k$.
\end{itemize}

In order to prove our main result, we need the following lemmas (cf. \cite{optimal}).

\begin{lemma}\label{sdfrew564}
Let  $\mathbb{E}\subseteq \mathbb{D}$ be two division rings  with $|\mathbb{E}|>2$, and let
$m,n, k\geq 2$ be integers with $m,n\geq k$.
Suppose $A\in\mathcal{M}_{k,k}(\mathbb{E}^{m\times n})$ and $(X,Y)$ is a fixed matrix representation of a point of $\mathscr{G}^l_{m+n-1,m-1}(\mathbb{D})$, where
$X\in \mathbb{D}^{m\times m}$ and $Y\in \mathbb{D}^{m\times n}$. Assume further that
${\rm ad}((X, Y), (I_m,B))=1$  for every $B\in\mathcal{M}_{k-1,k-1}(\mathbb{E}^{m\times n})$ satisfying ${\rm rank}(B-A)=1$, where $(I_m,B)$ is
a fixed matrix representation of some point of $\mathscr{G}^l_{m+n-1,m-1}(\mathbb{D})$.
Then $X$ is  invertible  and $Y=XA$.
In the case of $k=2$ there is  the additional possibility that $X$ is invertible  and $Y=0$.
\end{lemma}
\proof
{\bf Step 1.}
There is a permutation matrix $P$ and an invertible matrix $Q$ over $\mathbb{E}$,  such that $PAQ={\rm diag}(I_k,0)$.
Let $B\in\mathcal{M}_{k-1,k-1}(\mathbb{E}^{m\times n})$. Then $B\sim A$ if and only if $PBQ\sim PAQ$ with $PBQ\in\mathcal{M}_{k-1,k-1}(\mathbb{E}^{m\times n})$.
On the other hand, by (\ref{Grassmanndistance}),
 $\scriptsize {\rm ad}((X, Y),(I_{m}, B))=1$ if and only if ${\rm ad}((PXP^{-1}, PYQ),(I_{m}, PBQ))=1$. Moreover, $Y=XA$
 if and only if $PYQ=(PXP^{-1})(PAQ)$. Thus,  we may assume with no loss of generality that
$$A={\rm diag}(I_k,0).$$

Since ${\rm ad}((X, Y), (I_m,B))=1$ for any $B\in\mathcal{M}_{k-1,k-1}(\mathbb{E}^{m\times n})$ satisfying $B\sim A$,
from (\ref{Grassmanndistance}) we have
\begin{equation}\label{bcv36etertet}
 \mbox{ $Y-XB\sim 0$, \ for any $B\in\mathcal{M}_{k-1,k-1}(\mathbb{E}^{m\times n})$ satisfying $B\sim A$.}
\end{equation}

For any $B\in\mathcal{M}_{k-1,k-1}(\mathbb{E}^{m\times n})$ satisfying $B\sim A$, we get
 $m={\rm rank}(X,Y)={\rm rank}(X,Y-XB)\leq {\rm rank}(X)+{\rm rank}(Y-XB)$. It follows that ${\rm rank}(X)\geq m-1$.

We prove this lemma only for the case of $2\leq k<m$;  the case of $k=m$ is similar. From now on we assume that  $2\leq k<m$.
Write $X=(X_1,X_2)$ where $X_1\in\mathbb{D}^{m\times k}$ and  $X_2\in\mathbb{D}^{m\times (m-k)}$. Then $X_1\neq 0$ because $k\geq 2$.
Let ${\rm rank}(X_1)=s$ where $1\leq s\leq k$. Then there is  $P_1\in GL_m(\mathbb{E})$ such that $\scriptsize X_1=P_1\left(
                             \begin{array}{c}
                               X_{11} \\
                               0 \\
                             \end{array}
                           \right)$ where $X_{11}\in\mathbb{E}_s^{s\times k}$.
Without loss of generality,  we assume  $P_1=I_m$. Then $\scriptsize X=\left(
                             \begin{array}{cc}
                               X_{11}&X_{12} \\
                               0&X_{22} \\
                             \end{array}
                           \right)$ where $\scriptsize \left(
                             \begin{array}{c}
                               X_{12} \\
                               X_{22} \\
                             \end{array}
                           \right)=X_2$.
For any $B\in\mathcal{M}_{k-1,k-1}(\mathbb{E}^{m\times n})$ satisfying $B\sim A={\rm diag}(I_k,0)$, we have
$\scriptsize B=\left(
                             \begin{array}{c}
                               C \\
                               0 \\
                             \end{array}
                           \right)$, where $C\in\mathcal{M}_{k-1,k-1}(\mathbb{E}^{k\times n})$ satisfying $C\sim (I_k,0)$.
Write $\scriptsize Y=\left(
                             \begin{array}{c}
                               Y_1 \\
                               Y_2 \\
                             \end{array}
                           \right)$ where $Y_1\in \mathbb{D}^{s\times n}$ and $Y_2\in \mathbb{D}^{(m-s)\times n}$.
Then (\ref{bcv36etertet}) implies that
\begin{equation}\label{b76retbcvmv9i}
 \mbox{ $Y=\left(
  \begin{array}{c}
  Y_1\\
   Y_2 \\
   \end{array}
   \right)\sim \left(
  \begin{array}{c}
  X_{11}C \\
   0 \\
   \end{array}
   \right)$,  for any $C\in\mathcal{M}_{k-1,k-1}(\mathbb{E}^{k\times n})$ satisfying $C\sim (I_k,0)$.}
\end{equation}

For any $Z, W\in \mathbb{E}^{(k-1)\times 1}$,  let
$$ C_Z=\left(%
\begin{array}{ccc}
I_{k-1}&Z &0 \\
0&0&0\\
\end{array}%
\right), \ \ C_W'=\left(%
\begin{array}{ccc}
0&0&0\\
W&I_{k-1}&0 \\
\end{array}%
\right)\in \mathcal{M}_{k-1,k-1}(\mathbb{E}^{k\times n}).$$
Write $X_{11}=(x_{ij})_{s\times k}$, $Z=\,^t(z_1,\ldots,z_{k-1})$, $W=\,^t(w_2,\ldots,w_{k})$. Then we have
$$X_{11}C_Z=\left(\begin{array}{cccccccc}
x_{11} &x_{12} & \ldots & x_{1,k-1} & \sum_{j=1}^{k-1}x_{1j}z_j & 0 & \cdots &0 \\
x_{21} & x_{22} &\ldots & x_{2,k-1} & \sum_{j=1}^{k-1}x_{2j}z_j & 0 & \cdots &0 \\
 \vdots & \vdots& & \vdots & \vdots & \vdots &  &\vdots \\
 x_{s1} &x_{s2} & \ldots & x_{s,k-1} & \sum_{j=1}^{k-1}x_{sj}z_j & 0 & \cdots &0 \\
              \end{array}
            \right),$$
$$X_{11}C_W'=\left(\begin{array}{cccccccc}
 \sum_{j=2}^{k}x_{1j}w_j  & x_{12}&x_{13} &\ldots& x_{1k}& 0 & \cdots &0 \\
 \sum_{j=2}^{k}x_{2j}w_j & x_{22} &x_{23} &\ldots& x_{2k}& 0 & \cdots &0 \\
 \vdots & \vdots & \vdots &  &\vdots & \vdots &  &\vdots \\
  \sum_{j=2}^{k}x_{sj}w_j & x_{s2} &x_{s3} &\ldots& x_{sk}& 0 & \cdots &0 \\
              \end{array}
            \right).$$
 Since $C_Z\sim (I_k,0)$ and $C_W'\sim (I_k,0)$,
 from (\ref{b76retbcvmv9i}) we obtain that
\begin{equation}\label{53gdkltbcvmi}
 \mbox{ $Y=\left(
  \begin{array}{c}
  Y_1\\
   Y_2 \\
   \end{array}
   \right)\sim \left(
  \begin{array}{c}
  X_{11}C_Z \\
   0 \\
   \end{array}
   \right)$,  \  $Y=\left(
  \begin{array}{c}
  Y_1\\
   Y_2 \\
   \end{array}
   \right)\sim \left(
  \begin{array}{c}
  X_{11}C_W' \\
   0 \\
   \end{array}
   \right)$,  for all $Z, W\in \mathbb{E}^{(k-1)\times 1}$.}
\end{equation}

 {\bf Step 2.}  In this step,
we first show that $X_{11}C_Z\neq 0$ and $X_{11}C_W'\neq 0$ for all $Z, W\in \mathbb{E}^{(k-1)\times 1}$ by contradiction.
Suppose that $X_{11}C_Z=0$. Then, by $X_{11}\neq 0$ we get
$$X_{11}C_W'=\left(\begin{array}{cccccccc}
 x_{1k}w_k  &0&\ldots& 0&x_{1k}& 0 & \cdots &0 \\
 x_{2k}w_k & 0 &\cdots &0& x_{2k} &0 & \cdots &0 \\
 \vdots & \vdots &  & \vdots &\vdots & \vdots &  &\vdots \\
 x_{sk}w_k &0  &\ldots& 0& x_{sk} &0& \cdots &0 \\
              \end{array}
            \right)\neq 0. $$
Note that $Y\sim 0$ because (\ref{53gdkltbcvmi}).
Taking distinct $w_k$, by  (\ref{53gdkltbcvmi}) and  Lemma \ref{Matrix-PID5-4bb}, it is easy to see  that  $Y_2=0$, and hence ${\rm rank}(X_1,Y)={\rm rank}(X_{11}, Y_1)=s$.
Thus $m={\rm rank}(X,Y)\leq s+{\rm rank}(X_2)\leq s+m-k$, which implies that $k\leq s$. Since $k\geq s$, we obtain $k=s$.
Hence $X_{11}$ is invertible, a contradiction to $X_{11}C_Z=0$.
Therefore, $X_{11}C_Z\neq 0$ for all $Z\in \mathbb{E}^{(k-1)\times 1}$. Similarly,  $X_{11}C_W'\neq 0$ for all $W\in \mathbb{E}^{(k-1)\times 1}$.

Next, we will prove that $Y_2=0$ and $X$ is invertible.

Clearly, there are distinct $Z_i$ and $W_i\in\mathbb{E}^{(k-1)\times 1}$, $i=1,2$, such that $X_{11}C_{Z_1}\neq X_{11}C_{Z_2}$ and $X_{11}C_{W_1}'\neq X_{11}C_{W_2}'$.
Applying (\ref{53gdkltbcvmi}) and  Lemma \ref{Matrix-PID5-4bb}, it is easy to see  that $Y_2=0$. Thus ${\rm rank}(X_1,Y)={\rm rank}(X_{11}, Y_1)=s$.
By $m={\rm rank}(X,Y)\leq s+{\rm rank}(X_2)\leq s+m-k$, we get $k\leq s$ and hence  $k=s$.  By $k=s$ and $Y_2=0$, it is clear that $\scriptsize X=\left(
  \begin{array}{cc}
  X_{11} &X_{21}\\
   0&X_{22} \\
   \end{array}
   \right)$ is invertible.

Finally, we complete our proof.
Without loss of generality, we may assume  $X_{11}=I_k$. Thus, (\ref{53gdkltbcvmi}) and $Y_2=0$ imply that
\begin{equation}\label{cvgd53u57mi65}
 \mbox{$Y_1\sim C_Z$, \ \ $Y_1\sim C_W'$, \  for all $Z, W\in \mathbb{E}^{(k-1)\times 1}$.}
\end{equation}
When $k\geq 3$, (\ref{cvgd53u57mi65}) and Lemma \ref{Matrix-PID5-4bb} imply that $Y={\rm diag}(I_k,0)=A$, and hence $Y=XA=A$.
When $k=2$, applying (\ref{cvgd53u57mi65}) and Lemma \ref{Matrix-PID5-4bb}, we can obtain either $Y={\rm diag}(I_2,0)=A=XA$ or $Y=0$.
$\qed$

\begin{lemma}\label{sdfrxncdsdv8}  Let  $\mathbb{E}\subseteq \mathbb{D}$ be two division rings  with $|\mathbb{E}|>2$, and let
$m,n, k,r$ be integers with $m\geq k>r\geq 1$, $n> r$.
Suppose $A\in\mathcal{M}_{k,r}(\mathbb{E}^{m\times n})$ and $(X,Y)$ is a fixed matrix representation of a point of $\mathscr{G}^l_{m+n-1,m-1}(\mathbb{D})$, where
$X\in \mathbb{D}^{m\times m}$ and $Y\in \mathbb{D}^{m\times n}$. Assume further that
${\rm ad}((X, Y), (I_m,B))=1$  for every $B\in\mathcal{M}_{k,r+1}(\mathbb{E}^{m\times n})$ satisfying ${\rm rank}(B-A)=1$, where $(I_m,B)$ is
a fixed matrix representation of some point of $\mathscr{G}^l_{m+n-1,m-1}(\mathbb{D})$.
Then $X$ is invertible and $Y=XA$.
\end{lemma}
\proof{\bf Step 1.}  \
Using the same argument as in the proof of Lemma \ref{sdfrew564}, we may assume with no
  loss of generality that
$$A=\left(%
\begin{array}{cc}
 I_r &0\\
A_2&0\\
0&0\\
\end{array}%
\right)\in\mathcal{M}_{k,r}(\mathbb{E}^{k\times n}),$$
 where  $A_2\in\mathbb{E}^{(k-r)\times r}$ has no zero row.

Since $\scriptsize {\rm ad}((X, Y), (I_m,B))={\rm rank}\left(
 \begin{array}{cc}
 X & Y \\
 I_{m} & B \\
 \end{array}
 \right)-m=1$ for every $B\in\mathcal{M}_{k, r+1}(\mathbb{E}^{m\times n})$ with $B\sim A$,
we get
\begin{equation}\label{bceteg45r5353}
 \mbox{ $Y-XB\sim 0$, \ for any $B\in\mathcal{M}_{k,r+1}(\mathbb{E}^{m\times n})$ with $B\sim A$.}
\end{equation}
For any $B\in\mathcal{M}_{k,r+1}(\mathbb{E}^{m\times n})$ satisfying $B\sim A$, we have that
 $m={\rm rank}(X,Y)={\rm rank}(X,Y-XB)\leq {\rm rank}(X)+{\rm rank}(Y-XB)$. Thus ${\rm rank}(X)\geq m-1$.

We prove this lemma only for the case of $2\leq k<m$;  the case of $k=m$ is similar. From now on we assume that  $2\leq k<m$.
Let $X=(X_1,X_2)$, where $X_1\in\mathbb{D}^{m\times k}$, $X_2\in\mathbb{D}^{m\times (m-k)}$. Then $X_1\neq 0$ because $k\geq 2$.
Let ${\rm rank}(X_1)=s$ where $1\leq s\leq k$. Then, there is  $P_1\in GL_m(\mathbb{D})$ such that $\scriptsize X_1=P_1\left(
                             \begin{array}{cc}
                               X_{11}&X_{12}\\
                               0 &0 \\
                             \end{array}
                           \right)$ where $X_{11}\in\mathbb{D}^{s\times r}$ and $X_{12}\in\mathbb{D}^{s\times (k-r)}$.
Without loss of generality,  we assume that $P_1=I_m$. Then
$$X=(X_1,X_2)=\left(\begin{array}{ccc}
                               X_{11}&X_{12}&X_{21}\\
                               0 &0 &X_{22} \\
                             \end{array}
                           \right)$$
 and ${\rm rank}(X)\leq s+m-k$, where $X_{22}\in \mathbb{D}^{(m-s)\times(m-k)}$.
Write $\scriptsize Y=\left(
\begin{array}{c}
  Y_1 \\
   Y_2 \\
   \end{array}
  \right)$ where $Y_1\in \mathbb{D}^{s\times n}$ and $Y_2\in \mathbb{D}^{(m-s)\times n}$.

 Put
$\small C_Z=\left(%
\begin{array}{cc}
 I_r &0\\
A_2&Z\\
0&0\\
\end{array}%
\right)$ where $Z\in \mathbb{E}_1^{(k-r)\times (n-r)}$. Then $C_Z\in \mathcal{M}_{k,r+1}(\mathbb{E}^{m\times n})$ and $C_Z\sim A$.
It follows from (\ref{bceteg45r5353}) that
\begin{equation}\label{bcadaqw13lwqw}
\mbox{$Y=\left(
  \begin{array}{c}
  Y_1\\
   Y_2 \\
   \end{array}
   \right)\sim \left(
  \begin{array}{cc}
  X_{11}+X_{12}A_2 & X_{12}Z\\
   0 &0 \\
   \end{array}
   \right)$, \  for all $Z\in \mathbb{E}_1^{(k-r)\times (n-r)}$.}
\end{equation}

Since $A_2=(a_{ij})$ has no zero row, we may assume with no loss of generality that $a_{11}\neq 0$.
Then $(a_{11}, \ldots, a_{1r})(I_r-\lambda E_{11}^{r\times r})^{-1}E_{11}^{r\times (n-r)}\neq 0$  for any $\lambda\in \mathbb{E}$ with $\lambda\neq 1$.
It follows that $A_2(I_r-\lambda E_{11}^{r\times r})^{-1}E_{11}^{r\times (n-r)}\neq 0$ for any $\lambda\in \mathbb{E}$ with $\lambda\neq 1$.
Thus, for any $\lambda\in \mathbb{E}$ with $\lambda\neq 1$, we have
$$E_{\lambda}:=\left(
              \begin{array}{cc}
                I_r-\lambda E_{11}^{r\times r} & E_{11}^{r\times (n-r)} \\
                A_2 & 0 \\
                0&0\\
              \end{array}
            \right)\in \mathcal{M}_{k,r+1}(\mathbb{E}^{m\times n})$$
and $E_{\lambda}\sim A$. By  (\ref{bceteg45r5353}), we get
\begin{equation}\label{31yryfhi6cv13lw}
\mbox{$\left(
  \begin{array}{c}
  Y_1\\
   Y_2 \\
   \end{array}
   \right)\sim \left(
  \begin{array}{cc}
  X_{11}(I_r-\lambda E_{11}^{r\times r})+X_{12}A_2 & X_{11}E_{11}^{r\times (n-r)}\\
   0 &0 \\
   \end{array}
   \right)$,  for any $\lambda\in \mathbb{E}$ with $\lambda\neq 1$.}
\end{equation}

For any $\lambda\in \mathbb{E}$, let
$$C_\lambda=\left(
              \begin{array}{cc}
                I_r & 0 \\
                A_2 -\lambda E_{11}^{(k-r)\times r}&  E_{11}^{(k-r)\times (n-r)} \\
                0&0\\
              \end{array}
            \right)\in \mathcal{M}_{k,r+1}(\mathbb{E}^{m\times n}).$$
Then $C_\lambda\sim A$. Hence (\ref{bceteg45r5353}) implies that
\begin{equation}\label{vc18679g0fs46}
\mbox{$Y=\left(
  \begin{array}{c}
  Y_1\\
   Y_2 \\
   \end{array}
   \right)\sim \left(
  \begin{array}{cc}
  X_{11}+X_{12}( A_2 -\lambda E_{11}^{(k-r)\times r}) & X_{12}E_{11}^{(k-r)\times (n-r)}\\
   0 &0 \\
   \end{array}
   \right)$,  for any $\lambda\in \mathbb{E}$.}
\end{equation}

{\bf Step 2.} \ In this step,  we first show $X_{12}\neq 0$ by contradiction.
Suppose  $X_{12}=0$. Then  $s\leq r$ and (\ref{bcadaqw13lwqw}) implies that $Y-XA\sim 0$.
Thus $m={\rm rank}(X,Y)={\rm rank}(X,Y-XA)\leq {\rm rank}(X)+{\rm rank}(Y-XA)\leq s+m-k+1$, and hence $r\geq s\geq k-1\geq r$.
Consequently $s=r=k-1$ and $X_{11}$ is invertible. Without loss of generality,  we may assume that $X_{11}=I_r$.
Then $\small X=(X_1,X_2)=\left(
                           \begin{array}{ccc}
                             I_r & 0 &X_{21}\\
                             0 & 0 &X_{22}\\
                           \end{array}
                         \right)$ where $X_{22}\in \mathbb{D}^{(m-r)\times(m-k)}$.
By (\ref{31yryfhi6cv13lw}) we have
\begin{equation}\label{jgasmvb74egsg}
\mbox{ $Y=\left(
  \begin{array}{c}
  Y_1\\
   Y_2 \\
   \end{array}
   \right)\sim \left(
  \begin{array}{cc}
  I_r-\lambda E_{11}^{r\times r} &  E_{11}^{r\times (n-r)}\\
   0 &0 \\
   \end{array}
   \right)$, \  for any $\lambda\in \mathbb{E}$ with $\lambda\neq 1$.}
\end{equation}
From (\ref{bcadaqw13lwqw}) we get $Y\sim {\rm diag}(I_r, 0)$. Thus,
using  Lemma \ref{Matrix-PID5-4bb} and (\ref{jgasmvb74egsg}), we can prove that $Y_2=0$. Then ${\rm rank}(X_1,Y)=r$.
By $m={\rm rank}(X,Y)\leq r+{\rm rank}(X_{22})\leq r+m-k$, we get $k\leq r$, a contradiction to the condition $k>r$.
Therefore, we have proved  $X_{12}\neq 0$.

Next, we affirm that $Y_2=0$ and $X$ is invertible.
Clearly,  there are two distinct $Z_1,Z_2\in \mathbb{E}_1^{(k-r)\times (n-r)}$ such that $X_{12}Z_1\neq X_{12}Z_2$.
By (\ref{bcadaqw13lwqw}), we obtain
$$ Y=\left(
  \begin{array}{c}
  Y_1\\
   Y_2 \\
   \end{array}
   \right)\sim \left(
  \begin{array}{cc}
  X_{11}+X_{12}A_2 & X_{12}Z_i\\
   0 &0 \\
   \end{array}
   \right), \ i=1,2.$$
 Thus Lemma \ref{Matrix-PID5-4bb} implies that $Y_2=0$, or $\small Y=\left(
  \begin{array}{cc}
  X_{11}+X_{12}A_2 &Y_{12}\\
   0&Y_{22} \\
   \end{array}
   \right)$ where $Y_{22}\in\mathbb{D}^{(m-s)\times (n-r)}$, $(0,Y_{22})=Y_2$ and $(X_{11}+X_{12}A_2, Y_{12})=Y_1$.

On the other hand,
we have  $X_{11}\lambda E_{11}^{r\times r}\neq 0$ for any $\lambda\in\mathbb{E}^*$, or $X_{12}\lambda E_{11}^{(k-r)\times r}\neq 0$ for any $\lambda\in\mathbb{E}^*$.
 Otherwise, the first column of $X_{11}$ and the first column of $X_{12}$ are zeros, thus
$m-1\leq {\rm rank}(X)\leq m-2$, a contradiction. By (\ref{31yryfhi6cv13lw}) and (\ref{vc18679g0fs46}), we may assume with no loss of generality that
 $X_{11}\lambda E_{11}^{r\times r}\neq 0$ for any $\lambda\in\mathbb{E}^*$.
Suppose   $\small Y=\left(
  \begin{array}{cc}
  X_{11}+X_{12}A_2 &Y_{12}\\
   0&Y_{22} \\
   \end{array}
   \right)$.  From (\ref{31yryfhi6cv13lw}) we get
$$\left(
  \begin{array}{cc}
  X_{11}+X_{12}A_2 &Y_{12}\\
   0&Y_{22} \\
   \end{array}
   \right)\sim \left(
  \begin{array}{cc}
  X_{11}(I_r-\lambda E_{11}^{r\times r})+X_{12}A_2 & X_{11}E_{11}^{r\times (n-r)}\\
   0 &0 \\
   \end{array}
   \right),$$
for any $\lambda\in\mathbb{E}^*$ with $\lambda\neq 1$. Since $ X_{11}+X_{12}A_2 \neq X_{11}(I_r-\lambda E_{11}^{r\times r})+X_{12}A_2$ for any $\lambda\in\mathbb{E}^*$ with $\lambda\neq 1$,
we have $Y_{22}=0$. Thus, we always have $Y_2=0$. By $\scriptsize (X, Y)=\left(
                           \begin{array}{cccc}
                             X_{11} & X_{12} &X_{21}& Y_1\\
                             0 & 0 &X_{22}&0\\
                           \end{array}
                         \right)$, it is clear that $X$ is invertible.

By $m={\rm rank}(X)\leq s+m-k$, we get $k\leq s$. Since $k\geq s$, we obtain $k=s$. Thus $X_{11}\in\mathbb{D}^{k\times r}$, $X_{12}\in\mathbb{D}^{k\times (k-r)}$
and $X_{22}\in \mathbb{D}^{(m-k)\times(m-k)}$. Without loss of generality, we assume  $X=I_m$. Then $\small X_{11}=\left(
                                                                                                                           \begin{array}{c}
                                                                                                                             I_r \\
                                                                                                                             0 \\
                                                                                                                           \end{array}
                                                                                                                         \right)$,
                                                                                                                          $\small X_{12}=\left(
                                                                                                                           \begin{array}{c}
                                                                                                                             0 \\
                                                                                                                             I_{k-r} \\
                                                                                                                           \end{array}
                                                                                                                         \right)$,
$X_{21}=0$ and $X_{22}=I_{m-k}$.

Finally, we assert $Y=XA$.
By (\ref{bcadaqw13lwqw}) we have $\scriptsize Y-A\sim
\left(\begin{array}{cc}
 0_r &0\\
 0 &Z \\
 0 &0\\
 \end{array}
 \right)$  for all $Z\in\mathbb{E}_1^{(k-r)\times (n-r)}$. By  Lemma \ref{Matrix-PID5-4bb}, we get
 $\scriptsize Y-A=\left(\begin{array}{cc}
 0_r &W_1\\
 0 &W_3 \\
  0 &W_4 \\
 \end{array}
 \right)$ or $\scriptsize Y-A=\left(\begin{array}{cc}
 0_r & 0\\
 W_2 &W_3 \\
  0 &0 \\
 \end{array}
 \right)$, where $W_3\in\mathbb{D}^{(k-r)\times (n-r)}$. By (\ref{31yryfhi6cv13lw}), we have
 $$\mbox{$ Y-A\sim
\left(\begin{array}{cc}
 -\lambda E_{11}^{r\times r} &E_{11}^{r\times (n-r)}\\
   0 &0 \\
    0 &0 \\
 \end{array}
 \right)$, \  for any $\lambda\in \mathbb{E}$ with $\lambda\neq 1$.}$$
  It follows that
 $\scriptsize Y-A=\left(\begin{array}{cc}
 0_r &W_1\\
 0 &0 \\
  0 &0 \\
 \end{array}
 \right)$. From (\ref{vc18679g0fs46}) we obtain  that
 $$Y-A=\left(\begin{array}{cc}
 0_r &W_1\\
 0 &0 \\
  0 &0 \\
 \end{array}
 \right)\sim \left(\begin{array}{cc}
 0_r &0\\
 -\lambda E_{11}^{(k-r)\times r}  &E_{11}^{(k-r)\times (n-r)} \\
 0 &0\\
 \end{array}
 \right),$$
 for any $\lambda\in \mathbb{E}$, and hence $W_1=0$. Then, we have proved $Y=A=XA$.
$\qed$

Similarly, we can prove the following lemmas (cf. \cite{optimal}).

\begin{lemma}\label{sd5523w564}
 Let  $\mathbb{E}\subseteq \mathbb{D}$ be two division rings  with $|\mathbb{E}|>2$, and let $m,n, k\geq 2$ be integers with $m, n\geq k$.
 Suppose $A\in\mathcal{N}_{k,k}(\mathbb{E}^{m\times n})$ and
$\scriptsize\left(
\begin{array}{c}
X \\
Y\\
\end{array}
\right)$ is a fixed matrix representation of a point of $\mathscr{G}^r_{m+n-1,n-1}(\mathbb{D})$, where
$X\in \mathbb{D}^{n\times n}$ and $Y\in \mathbb{D}^{m\times n}$. Assume further that
$\scriptsize {\rm ad}\left(\left(
\begin{array}{c}
X \\
Y\\
\end{array}
\right), \left(
\begin{array}{c}
I_{n} \\
 B\\
 \end{array}
 \right)\right)=1$  for every $B\in\mathcal{N}_{k-1,k-1}(\mathbb{E}^{m\times n})$ satisfying ${\rm rank}(B-A)=1$, where $\scriptsize\left(
\begin{array}{c}
I_{n} \\
 B\\
 \end{array}
 \right)$ is
a fixed matrix representation of some point of $\mathscr{G}^r_{m+n-1,n-1}(\mathbb{D})$.
Then $X$ is  invertible  and $Y=AX$.
In the case of $k=2$ there is  the additional possibility that $X$ is invertible  and $Y=0$.
\end{lemma}

\begin{lemma}\label{sdfr533sdv8}
Let  $\mathbb{E}\subseteq \mathbb{D}$ be two division rings  with $|\mathbb{E}|>2$, and let
 $m,n, k,r$ be integers with $n\geq k>r\geq 1$, $m> r$.
Suppose $A\in\mathcal{N}_{k,r}(\mathbb{E}^{m\times n})$ and $\scriptsize\left(
\begin{array}{c}
X \\
Y\\
\end{array}
\right)$ is a fixed matrix representation of a point of $\mathscr{G}^r_{m+n-1,n-1}(\mathbb{D})$, where
$X\in \mathbb{D}^{n\times n}$ and $Y\in \mathbb{D}^{m\times n}$. Assume further that
$\scriptsize {\rm ad}\left(\left(
\begin{array}{c}
X \\
Y\\
\end{array}
\right), \left(
\begin{array}{c}
I_{n} \\
 B\\
 \end{array}
 \right)\right)=1$   for every $B\in\mathcal{N}_{k,r+1}(\mathbb{E}^{m\times n})$ satisfying ${\rm rank}(B-A)=1$, where $\scriptsize\left(
\begin{array}{c}
I_{n} \\
 B\\
 \end{array}
 \right)$  is a fixed matrix representation of some point of $\mathscr{G}^l_{m+n-1,n-1}(\mathbb{D})$.
Then $X$ is invertible and $Y=AX$.
\end{lemma}

\section{Proofs of Theorem \ref{MainTheorem0} and Corollary \ref{MainTheoremcor0}}

In this section, we will prove Theorem \ref{MainTheorem0} and Corollary \ref{MainTheoremcor0}.

\begin{lemma}\label{Grassmann-properties3}   Suppose  $\mathbb{D}, \mathbb{D}'$ are division rings  and there exists a nonzero ring homomorphism $\tau$
from $\mathbb{D}$ to $\mathbb{D}'$. Let $m,n\geq 2$ be integers and $L\in {\mathbb{D}'}^{n\times m}$.
Then $I_m+X^\tau L\in GL_m(\mathbb{D}')$ $(X\in\mathbb{D}^{m\times n})$ if and only if
$I_n+LX^\tau\in GL_n(\mathbb{D}')$ $(X\in\mathbb{D}^{m\times n})$. Moreover, if $I_m+X^\tau L\in GL_m(\mathbb{D}')$ $(X\in\mathbb{D}^{m\times n})$, then
\begin{equation}\label{Grassmann-p-03}
\mbox{$(I_m+X^\tau L)^{-1}X^\tau=X^\tau(I_n+LX^\tau)^{-1}$, \ $X\in\mathbb{D}^{m\times n}$.}
\end{equation}
\end{lemma}
\proof    We have
\begin{equation}\label{Grassmann-p-03a}
\mbox{$(I_m+X^\tau L)X^\tau=X^\tau(I_n+LX^\tau)$, \ for all $X\in\mathbb{D}^{m\times n}$.}
\end{equation}
Let  $G_X=I_m+X^\tau L$, $X\in\mathbb{D}^{m\times n}$.
Assume that $G_X\in GL_m(\mathbb{D}')$ for all $X\in\mathbb{D}^{m\times n}$. By (\ref{Grassmann-p-03a}), we get
 $X^\tau=G_X^{-1}X^\tau(I_n+LX^\tau)$. Hence $I_n=(I_n-LG_X^{-1}X^\tau)(I_n+LX^\tau)$. Consequently, $I_n+LX^\tau\in GL_n(\mathbb{D}')$
 for all $X\in\mathbb{D}^{m\times n}$.
Similarly,  if $I_n+LX^\tau\in GL_n(\mathbb{D}')$ for all $X\in\mathbb{D}^{m\times n}$,  then $I_m+X^\tau L\in GL_m(\mathbb{D}')$ for all $X\in\mathbb{D}^{m\times n}$.
Thus,  $I_m+X^\tau L\in GL_m(\mathbb{D}')$  ($X\in\mathbb{D}^{m\times n}$) if and only if
$I_n+LX^\tau\in GL_n(\mathbb{D}')$ ($X\in\mathbb{D}^{m\times n}$). Moreover, if $I_m+X^\tau L\in GL_m(\mathbb{D}')$ for all $X\in\mathbb{D}^{m\times n}$,
then (\ref{Grassmann-p-03a}) implies that (\ref{Grassmann-p-03}) holds.
$\qed$

Similarly, we have the following lemma.

\begin{lemma}\label{Grassmann-properties4}   Suppose that $\mathbb{D}, \mathbb{D}'$ are division rings and there exists
a nonzero ring anti-homomorphism $\sigma$ from $\mathbb{D}$ to $\mathbb{D}'$. Let $m,n\geq 2$ be
integers and  $L\in {\mathbb{D}'}^{m\times n}$.
Then   $I_m+L\,^tX^\sigma\in GL_m(\mathbb{D}')$ $(X\in\mathbb{D}^{m\times n})$ if and only if
$I_n+ \,^tX^\sigma L\in GL_n(\mathbb{D}')$ $(X\in\mathbb{D}^{m\times n})$. Moreover, if $I_m+L\,^tX^\sigma\in GL_m(\mathbb{D}')$ $(X\in\mathbb{D}^{m\times n})$, then
$$
\mbox{$^tX^\sigma(I_m+L\,^tX^\sigma)^{-1}=(I_n+ \,^tX^\sigma L)^{-1}\,^tX^\sigma$, \ $X\in\mathbb{D}^{m\times n}$.}
$$
\end{lemma}

\begin{proposition}\label{MainLemma001}   Let $\mathbb{D}, \mathbb{D}'$ be division rings with $|\mathbb{D}|\geq 4$,  and let $m,n,m',n'\geq 2$ be
integers with $m'\geq m$, $n'\geq n$. Suppose that  $\theta:  {\mathbb{D}}^{m\times n}\rightarrow  {\mathbb{D}'}^{m'\times n'}$ is
a non-degenerate graph homomorphism with $\theta(0)=0$. Assume further that $\theta(\mathcal{M}_i)\subseteq \mathcal{M}'_i$, $i=1,\ldots,m$,
and $\theta(\mathcal{N}_j)\subseteq \mathcal{N}'_j$, $j=1,\ldots,n$.
Then there exist two invertible diagonal matrices  $P\in GL_{m'}(\mathbb{D}')$ and  $Q\in GL_{n'}(\mathbb{D}')$,
a nonzero ring homomorphism $\tau: \mathbb{D} \rightarrow\mathbb{D}'$, and a matrix $L\in {\mathbb{D}'}^{n\times m}$ with the property
that $I_m+X^\tau L\in GL_m(\mathbb{D}')$ for every $X\in {\mathbb{D}}^{m\times n}$, such that
\begin{equation}\label{35hfhbcnnv7}
\theta(X)=P\left(\begin{array}{cc}
 (I_m+X^\tau L)^{-1}X^\tau  & 0 \\
  0 & 0 \\
   \end{array}
  \right)Q, \ \, X\in {\mathbb{D}}^{m\times n}.
\end{equation}
In particular, if $\tau$ is surjective,  then $L=0$ and  $\tau$ is a ring isomorphism.
\end{proposition}
\proof
We  prove this result only for the case $m\leq n$.
 When $m>n$, by the symmetry of rows and columns of a matrix, we can prove similarly
$$\theta(X)=P\left(\begin{array}{cc}
 X^\tau (I_n+LX^\tau)^{-1} & 0 \\
  0 & 0 \\
   \end{array}
  \right)Q, \ \, X\in {\mathbb{D}}^{m\times n},$$
 where $P\in GL_{m'}(\mathbb{D}')$ and  $Q\in GL_{n'}(\mathbb{D}')$
are two invertible diagonal matrices, $\tau: \mathbb{D} \rightarrow\mathbb{D}'$ is a nonzero ring homomorphism,
$L\in {\mathbb{D}'}^{n\times m}$ with the property
that $I_n+LX^\tau \in GL_n(\mathbb{D}')$ for every $X\in {\mathbb{D}}^{m\times n}$. By Lemma \ref{Grassmann-properties3} ,
we have (\ref{35hfhbcnnv7}). From now on we assume that $m\leq n$.

{\bf Step 1.} \ For $1\leq i \leq m$,
by the conditions and Lemma \ref{degenerate-2},  the restriction map $\theta\mid_{\mathcal{M}_i}: \mathcal{M}_i\rightarrow \mathcal{M}'_i$
is an injective  weighted semi-affine map.  By  Lemma \ref{weaksemilinea34ear2},
we have
\begin{equation}\label{trghj7we99}
\theta\left(
         \begin{array}{c}
         0_{i-1,n} \\
           x_i \\
           0 \\
         \end{array}
       \right)=
\left(\begin{array}{c}
0_{i-1,n} \\
(k_i(x_i))^{-1}x_i^{\sigma_i}Q_i \\
0 \\
\end{array}
\right),  \  x_i\in {\mathbb{D}}^n,  1\leq i \leq m,
\end{equation}
where  $\sigma_i: \mathbb{D}\rightarrow \mathbb{D}'$ is a nonzero ring homomorphism, $Q_i\in {\mathbb{D}'}^{n\times n'}$, and
$k_i(x_i)=\sum_{j=1}^nx_{ij}^{\sigma_i} a_{ji}+b_i\neq 0$ for all $x_i=(x_{i1}, \ldots, x_{in})\in {\mathbb{D}}^n$ with
 $a_{ji}, b_i\in \mathbb{D}'$ ($i=1,\ldots,m$, $j=1,\ldots,n$). Note that all $a_{ji}$, $b_i$ are fixed with $b_i\neq 0$, $i=1,\ldots,m$, $j=1,\ldots,n$.
  In (\ref{trghj7we99}), some zeros are absent when $i=1$ or $i=m$.

Let $E_{ij}=E_{ij}^{m\times n}$ and $E_{ij}'=E_{ij}^{m'\times n'}$. Note that $\mathcal{M}_i\cap \mathcal{N}_j=\mathbb{D}E_{ij}$.
Since $\theta(\mathcal{M}_i)\subseteq \mathcal{M}'_i$ and $\theta(\mathcal{N}_j)\subseteq \mathcal{N}'_j$, we get
\begin{equation}\label{vcx34fsf67h}
\mbox{$\theta(E_{ij})=e_{ij}E_{ij}'$ (where $e_{ij}\in {\mathbb{D}'}^*$), \ $i=1,\ldots,m$, $j=1,\ldots,n$.}
\end{equation}
Thus, it is easy to see that  $Q_{i}=(T_i, 0)$,  where $T_i={\rm diag}(q_{i1}, \ldots, q_{in})$ is an invertible diagonal matrix,
$i=1,\ldots,m$. Hence
\begin{equation}\label{trgh34fs42bmm}
\theta\left(
         \begin{array}{c}
         0_{i-1,n} \\
           x_i \\
           0 \\
         \end{array}
       \right)=
\left(\begin{array}{cc}
0_{i-1,n}&0 \\
(k_i(x_i))^{-1}x_i^{\sigma_i}T_i&0 \\
0 &0_{m'-i,n'-n}\\
\end{array}
\right),  \  x_i\in {\mathbb{D}}^n, 1\leq i \leq m.
\end{equation}
Replacing $\theta$ by the map $X\mapsto \theta(X){\rm diag}(T_1^{-1}, I_{n'-n})$, we have $T_1=I_n$.

For a nonzero $x_1=(x_{11}, \ldots, x_{1n})\in {\mathbb{D}}^n$ and an integer $i$ with $2\leq i\leq m$, since
$$\theta\left(
         \begin{array}{c}
           x_1 \\
           0 \\
           0\\
         \end{array}
       \right)\sim \theta\left(
         \begin{array}{c}
         0_{i-1,n} \\
           x_1 \\
           0 \\
         \end{array}
       \right),  $$
it follows from (\ref{trgh34fs42bmm}) that
$$\left(
         \begin{array}{c}
           (k_1(x_1))^{-1}x_1^{\sigma_1}\\
           0 \\
           0\\
         \end{array}
       \right)\sim \left(
         \begin{array}{c}
         0_{i-1,n} \\
           (k_i(x_1))^{-1}x_1^{\sigma_i}T_i \\
           0 \\
         \end{array}
       \right),  $$
and hence
\begin{equation}\label{vcx3dfhdfit54y}
\left(
         \begin{array}{c}
           x_1^{\sigma_1} \\
           0  \\
           0 \\
         \end{array}
       \right)\sim \left(
         \begin{array}{c}
         0_{i-1,n} \\
         x_1^{\sigma_i}T_i \\
           0 \\
         \end{array}
       \right).
\end{equation}
Thus, there exists a nonzero $c_i\in \mathbb{D}'$ such that $c_ix_1^{\sigma_1}=x_1^{\sigma_i}T_i$. Then
$$\mbox{$c_i(x_{11}^{\sigma_1},\ldots, x_{1n}^{\sigma_1})=(x_{11}^{\sigma_i}q_{i1},\ldots, x_{1n}^{\sigma_i}q_{in})$ \ for all $x_{1j}\in \mathbb{D}$}.$$
Taking $x_{11}=\cdots=x_{1n}=1$, we get $c_i=q_{i1}=\cdots=q_{in}=:q_i$.  Thus $x_1^{\sigma_1}=x_1^{\sigma_i}$ for all  $x_1\in {\mathbb{D}}^n$, $i=1,\ldots, m$.
Let $\tau=\sigma_1=\sigma_2=\cdots=\sigma_m$. Then $T_i=q_iI_n$ and $x^{\tau}T_i=x^{\tau}q_i$, $x\in {\mathbb{D}}^n$, $i=2,\ldots,m$.
Note that the map $x\mapsto q_i^{-1}x^{\tau}q_i=:x^{\tau'}$ ($x\in \mathbb{D}$) is a nonzero ring homomorphism
from $\mathbb{D}$ to $\mathbb{D}'$. Moreover,  $(k_i(x_i))^{-1}x_i^{\tau} T_i=(k_i(x_i))^{-1}x_i^{\tau}q_i
= (k_i'(x_i))^{-1}x_i^{\tau'}$, where $k_i'(x_i)=\sum_{j=1}^nx_{ij}^{\tau'} q_i^{-1}a_{ji}+q_i^{-1}b_i$.
Thus,
we may assume with no loss of of generality that  $T_i=I_n$ in (\ref{trgh34fs42bmm}), $i=1,\ldots, m$.
Then  (\ref{trgh34fs42bmm}) becomes
\begin{equation}\label{tr7wadsaaae9wq9}
\theta\left(
         \begin{array}{c}
         0_{i-1,n} \\
           x_i \\
           0 \\
         \end{array}
       \right)=
\left(\begin{array}{cc}
0_{i-1,n}&0 \\
(k_i(x_i))^{-1}x_i^{\tau}&0 \\
0 &0_{m'-i,n'-n}\\
\end{array}
\right),  \  x_i\in {\mathbb{D}}^n, 1\leq i \leq m.
\end{equation}

Note that $((k_i(x_i))^{-1}=b_i^{-1}(\sum_{i=1}^nx_{ij}^\tau a_{ji}b_i^{-1}+1)^{-1}$.
Modify the map $\theta$ by the map $X\longmapsto{\rm diag}(b_1, \ldots,b_m,1,\ldots,1)\theta(X)$.
 We can assume that
$b_i=1$ in $k_i(x_i)$ ($i=1, \ldots, m$). Thus
\begin{equation}\label{dhlttrtr3244}
k_i(x_i)=\sum_{j=1}^nx_{ij}^\tau a_{ji}+1\neq 0, \ \ \mbox{$x_i=(x_{i1},\ldots, x_{in})\in {\mathbb{D}}^n$, \ $i=1,\ldots,m$,}
\end{equation}
where all $a_{ji}$, $i=1,\ldots,m$, $j=1,\ldots,n$, are fixed.

{\bf Step 2.} \ In this step, we prove that
\begin{equation}\label{fsd32vm69gdg0}
\mbox{${\rm rank}(\theta(A))=2$, \   $A\in\mathcal{M}_{2,2}( {\mathbb{D}}^{m\times n})$.}
\end{equation}
Write $\alpha=\{1,2\}$. We  prove (\ref{fsd32vm69gdg0}) only for the case  $\mathcal{M}_{2,2}^\alpha( {\mathbb{D}}^{m\times n})$;
the other cases on $\mathcal{M}_{2,2}({\mathbb{D}}^{m\times n})$  is similar.
Now, we assume that
$\mathcal{M}_{2,2}( {\mathbb{D}}^{m\times n})=\mathcal{M}_{2,2}^\alpha( {\mathbb{D}}^{m\times n})$.

Let
$\scriptsize A=\left(
     \begin{array}{c}
       \alpha_1 \\
       \alpha_2 \\
       0 \\
     \end{array}
   \right)\in \mathcal{M}_{2,2}^\alpha( {\mathbb{D}}^{m\times n})$, where  $\alpha_1, \alpha_2\in \mathbb{D}^n$
are left linearly independent over $\mathbb{D}$.
Write
$\scriptsize \theta(A)=\left(
     \begin{array}{c}
       \alpha_1^* \\
       \alpha_2^* \\
       0^* \\
     \end{array}
   \right)$ where  $\alpha_1^*, \alpha_2^*\in{\mathbb{D}'}^{n'}$.
We distinguish the following two cases to prove (\ref{fsd32vm69gdg0}).

{\em Case 2.1}. \
$\scriptsize\left(
     \begin{array}{c}
       \alpha_2^* \\
       0^* \\
     \end{array}
   \right)\neq 0$. Without loss of generality, we  assume that $\alpha_2^*\neq 0$. Let $\alpha_\lambda=\alpha_1+\lambda\alpha_2$
   where $\lambda\in \mathbb{D}$. Since
 $\scriptsize \left(
     \begin{array}{c}
      \alpha_1+\lambda\alpha_2 \\
       0_{m-1,n} \\
     \end{array}
   \right) \sim A$ for all $\lambda\in \mathbb{D}$,  from (\ref{tr7wadsaaae9wq9}) we have
  $$\mbox{$\theta\left(
     \begin{array}{c}
      \alpha_1+\lambda\alpha_2 \\
       0 \\
       0 \\
     \end{array}
   \right)=\left(
     \begin{array}{cc}
      (k_1(\alpha_\lambda))^{-1}(\alpha_1^\tau+\lambda^\tau\alpha_2^\tau)&0_{1,n'-n} \\
       0 &0\\
       0 &0\\
     \end{array}
   \right) \sim\left(
     \begin{array}{c}
      \alpha_1^* \\
       \alpha_2^* \\
       0^* \\
     \end{array}
   \right)$, \  for all $\lambda\in \mathbb{D}.$}$$
 Thus
 $$\mbox{${\rm rank}\left(
     \begin{array}{c}
      (k_1(\alpha_\lambda))^{-1}(\alpha_1^\tau+\lambda^\tau\alpha_2^\tau, \,0_{1,n'-n})-\alpha_1^* \\
       -\alpha_2^* \\
        -0^* \\
     \end{array}
   \right)=1$, \   for all $\lambda\in \mathbb{D}.$}$$
It follows that $ (k_1(\alpha_\lambda))^{-1}(\alpha_1^\tau+\lambda^\tau\alpha_2^\tau, \, 0_{1,n'-n})=\alpha_1^* +\lambda^\rho\alpha_2^*$ for all $\lambda\in \mathbb{D}$,
where $\rho: \mathbb{D}\rightarrow \mathbb{D}'$ is a map. Since
$\theta\mid_{\mathcal{M}_1}: \mathcal{M}_1\rightarrow \mathcal{M}'_1$ is an injective map, $\rho$ is  injective.
There are two distinct $\lambda_1, \lambda_2\in \mathbb{D}$ such that $\alpha_1 +\lambda_1\alpha_2$
and $\alpha_1+\lambda_2\alpha_2$ are left linearly independent over $\mathbb{D}$. Clearly,
 $(k_1(\alpha_{\lambda_1}))^{-1}(\alpha_1^\tau+\lambda_1^\tau\alpha_2^\tau, \, 0_{1,n'-n})=\alpha_1^* +\lambda_1^{\rho}\alpha_2^*$ and
 $(k_1(\alpha_{\lambda_2}))^{-1}(\alpha_1^\tau+\lambda_2^\tau\alpha_2^\tau, \, 0_{1,n'-n})=\alpha_1^* +\lambda_2^{\rho}\alpha_2^*$
 are also left linearly independent over $\mathbb{D}'$.
 Consequently,  $\alpha_1^*, \alpha_2^*$ are  left linearly independent over $\mathbb{D}'$. Hence ${\rm rank}(\theta(A))=2$.

{\em Case 2.2}. \
$\scriptsize \left(
     \begin{array}{c}
       \alpha_2^* \\
       0^* \\
     \end{array}
   \right)= 0$. Then $\scriptsize\theta(A)=\left(
     \begin{array}{c}
       \alpha_1^* \\
       0 \\
       0 \\
     \end{array}
   \right)$. Let $\beta_\lambda=\alpha_2+\lambda\alpha_1$, $\lambda\in \mathbb{D}$. Since
 $\scriptsize \left(
     \begin{array}{c}
     0\\
      \alpha_2+\lambda\alpha_1 \\
       0_{m-2,n} \\
     \end{array}
   \right) \sim A$ for all $\lambda\in \mathbb{D}$,  it follows from (\ref{tr7wadsaaae9wq9}) that
$$\mbox{$\theta\left(
     \begin{array}{c}
     0\\
      \alpha_2+\lambda\alpha_1 \\
       0_{m-2,n} \\
     \end{array}
   \right)=\left(
     \begin{array}{c}
     0\\
     (k_2(\beta_\lambda))^{-1}(\alpha_2^\tau+\lambda^\tau\alpha_1^\tau, \,0_{1,n'-n})\\
       0_{m'-2,n'} \\
     \end{array}
   \right) \sim\left(
     \begin{array}{c}
      \alpha_1^* \\
       0 \\
       0_{m'-2,n'} \\
     \end{array}
   \right)$,}$$
for all $\lambda\in \mathbb{D}.$ Hence
\begin{equation}\label{vcx32hfh86ad}
\mbox{${\rm rank}\left(
     \begin{array}{c}
      \alpha_1^* \\
      (k_2(\beta_\lambda))^{-1}(\alpha_2^\tau+\lambda^\tau\alpha_1^\tau, \,0_{1,n'-n})\\
     \end{array}
   \right)=1$, \    for all $\lambda\in \mathbb{D}$.}
\end{equation}
Note that  $k_2(\beta_\lambda)=k_2(\alpha_2)+\lambda^\tau k_2(\alpha_1)-\lambda^\tau$ and
$\theta\mid_{\mathcal{M}_2}: \mathcal{M}_2\rightarrow \mathcal{M}'_2$ is an injective map. If $\lambda_1\neq \lambda_2$, then
$(k_2(\beta_{\lambda_1}))^{-1}(\alpha_2^\tau+\lambda_1^\tau\alpha_1^\tau)\neq (k_2(\beta_{\lambda_2}))^{-1}(\alpha_2^\tau+\lambda_2^\tau\alpha_1^\tau)$,
and hence $\lambda_1^\tau \neq \lambda_2^\tau$.
There are two distinct $\lambda_1, \lambda_2\in \mathbb{D}$ such that $\alpha_2 +\lambda_1\alpha_1$
and $\alpha_2+\lambda_2\alpha_1$ are left linearly independent over $\mathbb{D}$. It is easy to see that
 $(k_2(\beta_{\lambda_1}))^{-1}(\alpha_2^\tau+\lambda_1^\tau\alpha_1^\tau)$ and  $(k_2(\beta_{\lambda_2}))^{-1}(\alpha_2^\tau+\lambda_2^\tau\alpha_1^\tau)$
 are also left linearly independent over $\mathbb{D}'$.
 Thus, from (\ref{vcx32hfh86ad}) we must have $\alpha_1^* =0$, and hence $\theta(A)=0$.

Since $\scriptsize \left(
     \begin{array}{c}
       \alpha_2 \\
       \alpha_1 \\
       0 \\
     \end{array}
   \right)\sim A$,
$\scriptsize\theta\left(
     \begin{array}{c}
       \alpha_2 \\
       \alpha_1 \\
       0 \\
     \end{array}
   \right)
=:\left(
     \begin{array}{c}
       \beta_1 \\
       \beta_2 \\
       0' \\
     \end{array}
   \right)\sim \theta(A)=0$, where $\beta_1, \beta_2\in \mathbb{{D}'}^{n'}$. Then
$\scriptsize {\rm rank}\left(
     \begin{array}{c}
       \beta_1 \\
       \beta_2 \\
       0' \\
     \end{array}
   \right)=1$. By
$\scriptsize {\rm rank}\left(
     \begin{array}{c}
       \alpha_2 \\
       \alpha_1 \\
       0 \\
     \end{array}
   \right)=2$,
we can get $\scriptsize \left(
     \begin{array}{c}
       \beta_2 \\
       0' \\
     \end{array}
   \right)= 0$. Otherwise, if
$\scriptsize \left(
     \begin{array}{c}
       \beta_2 \\
       0' \\
     \end{array}
   \right)\neq 0$,  similar to the proof of Case 2.1, we have
$\scriptsize {\rm rank}\left(
     \begin{array}{c}
       \beta_1 \\
       \beta_2 \\
       0' \\
     \end{array}
   \right)=2$, a contradiction. Hence $\scriptsize \theta\left(
     \begin{array}{c}
       \alpha_2 \\
       \alpha_1 \\
       0 \\
     \end{array}
   \right)
=\left(
     \begin{array}{c}
       \beta_1 \\
       0 \\
       0 \\
     \end{array}
   \right)$. We have
 $\scriptsize \theta\left(
     \begin{array}{c}
    0\\
      \alpha_1+\lambda\alpha_2 \\
       0_{m-2,n} \\
     \end{array}
   \right)\sim \theta\left(
     \begin{array}{c}
       \alpha_2 \\
       \alpha_1 \\
       0 \\
     \end{array}
   \right)$ for all $\lambda\in \mathbb{D}$.
 Similar to the proof of $\theta(A)=0$ above, we can prove  $\scriptsize \theta\left(
     \begin{array}{c}
       \alpha_2 \\
       \alpha_1 \\
       0 \\
     \end{array}
   \right)=0$,  a contradiction. Therefore, Case 2.2 cannot occur.

Combining Case 2.1 with Case 2.2, we obtain (\ref{fsd32vm69gdg0}).

{\bf Step 3.} \ In this step, we prove that $\theta$ is of the form
$$\theta(X)={\rm diag}\left((I_m+X^\tau L)^{-1}X^\tau, \,0\right),$$
where $L\in {\mathbb{D}'}^{n\times m}$ is fixed with the property that $I_m+X^\tau L\in GL_m(\mathbb{D}')$ for any $X\in {\mathbb{D}}^{m\times n}$.

By (\ref{dhlttrtr3244}), let $L=(a_{ij})\in {\mathbb{D}'}^{n\times m}$. For any $x_i\in {\mathbb{D}}^n$,
 we have that
$$\mbox{$\small I_m+X^\tau L=\left(
      \begin{array}{ccc}
        I_{i-1} &0 & 0\\
        * & k_i(x_i)& * \\
         0 & 0 & I_{m-i} \\
      \end{array}
    \right)\in GL_m(\mathbb{D}')$,  if $\small X=\left(
         \begin{array}{c}
         0_{i-1,n} \\
           x_i \\
           0 \\
         \end{array}
       \right)\in \mathbb{D}^{m\times n}$.}
 $$
Thus (\ref{tr7wadsaaae9wq9}) can be written as
\begin{equation}\label{tr7w77ae9wq9}
\mbox{$\theta(X)={\rm diag}\left((I_m+X^\tau L)^{-1}X^\tau, \,0\right)$,   if $\scriptsize X=\left(
         \begin{array}{c}
         0_{i-1,n} \\
           x_i \\
           0 \\
         \end{array}
       \right)\in \mathbb{D}^{m\times n}$,  $1\leq i \leq m$.}
\end{equation}

Put
$$\mathscr{L}=\left\{X\in {\mathbb{D}}^{m\times n}: \, \mbox{$I_m+X^\tau L$ is invertible}\right\}.$$
Let $\mathcal{L}_0=\{0\}\subset \mathscr{L}$ and
$\mathcal{L}_k=\bigcup_{j=1}^k\mathcal{M}_{k,j}(\mathbb{D}^{m\times n})$, $k=1,\ldots, m.$
Clearly,  $\bigcup_{k=0}^m\mathcal{L}_k= {\mathbb{D}}^{m\times n}.$
Define the map $\psi: \mathscr{L}\rightarrow {\mathbb{D}'}^{m'\times n'}$ by
$$\mbox{$\psi(X)={\rm diag}\left((I_m+X^\tau L)^{-1}X^\tau, \,0\right)$, \ $X\in \mathscr{L}$.}$$
We will prove that $\mathscr{L}={\mathbb{D}}^{m\times n}$ and $\theta(X)=\psi(X)$ for all $X\in {\mathbb{D}}^{m\times n}$.

By (\ref{tr7w77ae9wq9}), $\mathcal{L}_1\subset \mathscr{L}$ and $\theta(X)=\psi(X)$
for any $X\in \mathcal{L}_1$. Suppose that $\mathcal{L}_{k-1}\subset \mathscr{L}$ and $\theta(X)=\psi(X)$
for any $X\in \mathcal{L}_{k-1}$ ($2\leq k \leq m$). Applying the mathematical induction,
we prove that $\mathcal{L}_{k}\subset \mathscr{L}$ and $\theta(X)=\psi(X)$
for any $X\in \mathcal{L}_{k}$ as follows.

Let $\mathbb{E}=\tau(\mathbb{D})$. Then $\mathbb{E}$ is a division subring of $\mathbb{D}'$  with $|\mathbb{E}|\geq 4$,
and the $\tau$ is a ring isomorphism from $\mathbb{D}$ to $\mathbb{E}$. Moreover, we have  $\tau(\mathbb{D}^{m\times n})=\mathbb{E}^{m\times n}$ and
 $\tau(\mathcal{M}_{k,j}(\mathbb{D}^{m\times n}))=\mathcal{M}_{k,j}(\mathbb{E}^{m\times n})$.

Let $A\in \mathcal{M}_{k,k}(\mathbb{D}^{m\times n})$  ($k\leq m$). Then $A^\tau\in \mathcal{M}_{k,k}(\mathbb{E}^{m\times n})$.
Clearly, $B\in \mathcal{M}_{k-1,k-1}(\mathbb{D}^{m\times n})$ with $B\sim A$ if and only if $B^\tau\in \mathcal{M}_{k-1,k-1}(\mathbb{E}^{m\times n})$
with $B^\tau\sim A^\tau$.
Let  $B\in \mathcal{M}_{k-1,k-1}(\mathbb{D}^{m\times n})$ with $B\sim  A$. Then $\theta(A)\sim  \theta(B)$ and $B^\tau\sim A^\tau$.
 By the induction hypothesis, we get $B\in \mathscr{L}$ and
 $$\theta(B)={\rm diag}\left((I_m+B^\tau L)^{-1}B^\tau, \,0\right).$$

 There are two matrices $B_1,B_2\in\mathcal{M}_{k-1,k-1}(\mathbb{D}^{m\times n})$ such that ${\rm rank}(B_1-B_2)\geq 2$ and $A\sim B_i$, $i=1,2$. By Example \ref{wett8345},
${\rm rank}(\theta(B_1)-\theta(B_2))\geq 2$, and hence
\begin{equation}\label{bvc524tetd32u}
{\rm rank}\left((I_m+B_1^\tau L)^{-1}B_1^\tau-(I_m+B_2^\tau L)^{-1}B_2^\tau \right)\geq 2.
\end{equation}
Let $\scriptsize\theta(A)=\left(
                  \begin{array}{cc}
                    A_1 & A_2 \\
                    A_3 & A_4 \\
                  \end{array}
                \right)$ where $A_1\in {\mathbb{D}'}^{m\times n}$.
 Since $\theta(A)\sim \theta(B_i)$, $i=1,2$, we get
 $${\rm rank}\left(
               \begin{array}{cc}
                 A_1-(I_m+B_i^\tau L)^{-1}B_i^\tau & A_2 \\
                 A_3 & A_4 \\
               \end{array}
             \right)=1, \ \ i=1,2.$$
 It follows from Lemma \ref{Matrix-PID5-4bb} that $A_4=0$ and either $A_2=0$ or $A_3=0$. Applying (\ref{bvc524tetd32u}), it is easy to see that
$A_2=0$ and $A_3=0$. Thus
$$\theta(A)=\left(
                  \begin{array}{cc}
                    A_1 & 0 \\
                    0 & 0 \\
                  \end{array}
                \right).$$

 For any
$B^\tau\in \mathcal{M}_{k-1,k-1}(\mathbb{E}^{m\times n})$ with $B^\tau\sim A^\tau$, we have $\theta(A)\sim \theta(B)$, and hence
$$((I_m+B^\tau L)^{-1}B^\tau, I_{m})\sim (A_1, I_{m}).$$
Matrices $((I_m+B^\tau L)^{-1}B^\tau, I_{m})$ and $(A_1, I_{m})$ are two matrix representations of two points of Grassmann space
$\mathcal{G}^l_{m+n-1,m-1}(\mathbb{D}')$. In the viewpoint of Grassmann space, from (\ref{Grassmanndistance}) we have that
${\rm ad}((I_m+B^\tau L)^{-1}B^\tau, I_{m}), (A_1, I_{m}))=1$
for any $B^\tau\in \mathcal{M}_{k-1,k-1}(\mathbb{E}^{m\times n})$ with $B^\tau\sim A^\tau$.
 Note that  $((I_m+B^\tau L)^{-1}B^\tau, I_{m})$ and $(B^\tau, I_m+B^\tau L)$ are two matrix representations of
the same point of the Grassmann space.  By (\ref{Grassmanndistance}), we get ${\rm ad}((B^\tau, I_m+B^\tau L), (A_1, I_{m}))=1$
for any $B^\tau\in \mathcal{M}_{k-1,k-1}(\mathbb{E}^{m\times n})$ with $B^\tau\sim A^\tau$.  Since
$$(B^\tau, I_{m}+B^\tau L)\left(
                        \begin{array}{cc}
                          -L & I_{n} \\
                          I_{m} & 0 \\
                        \end{array}
                      \right)=(I_{m}, B^\tau) \ \ {\rm and } \ \
(A_1, I_{m})\left(
                        \begin{array}{cc}
                          -L & I_{n} \\
                          I_{m} & 0 \\
                        \end{array}
                      \right)=(I_{m}-A_1L, \, A_1),$$
it follows from (\ref{Grassmanndistance})  that  ${\rm ad}((I_{m}-A_1L, \, A_1), (I_{m}, B^\tau))=1$
for any $B^\tau\in \mathcal{M}_{k-1,k-1}(\mathbb{E}^{m\times n})$ with $B^\tau\sim A^\tau$. Let $W=I_{m}-A_1L$.
It follows from Lemma \ref{sdfrew564}  that $W$ is invertible and  $A_1=WA^\tau$. Since
$$(A_1,I_m)=(A_1, W+A_1L)=(WA^\tau, W+WA^\tau L)=W(A^\tau, I_m+A^\tau L),$$
we have  $I_m=W(I_m+A^\tau L)$. Hence $W=(I_m+A^\tau L)^{-1}$ and $A_1=(I_m+A^\tau L)^{-1}A^\tau$. Then
$\theta(A)={\rm diag}\left((I_m+A^\tau L)^{-1}A^\tau, \, 0\right)$, $A\in \mathcal{M}_{k,k}(\mathbb{D}^{m\times n})$.
Therefore,
$$
\mbox{$\mathcal{M}_{k,k}(\mathbb{D}^{m\times n})\subset\mathscr{L}$ \ and \ $\theta(A)=\psi(A)$ \
for any $A\in \mathcal{M}_{k,k}(\mathbb{D}^{m\times n})$.}
$$

By Lemma \ref{sdfrxncdsdv8}, we have similarly that
$$
\mbox{$\mathcal{M}_{k,k-1}(\mathbb{D}^{m\times n})\subset\mathscr{L}$ \ and \ $\theta(X)=\psi(X)$ \
for any $X\in \mathcal{M}_{k,k-1}(\mathbb{D}^{m\times n})$.}
$$
Using this method, we can prove  $\mathcal{M}_{k,j}(\mathbb{D}^{m\times n})\subset\mathscr{L}$ and $\theta(A)=\psi(A)$
for any $A\in \mathcal{M}_{k,j}(\mathbb{D}^{m\times n})$,
$j=k-2,\ldots, 1$. Thus,  $\mathcal{L}_{k}\subset\mathscr{L}$ and $\theta(X)=\psi(X)$ for any $X\in \mathcal{L}_{k}$.
By the mathematical induction, we have  that $\mathscr{L}={\mathbb{D}}^{m\times n}$ and
$$\theta(X)=\psi(X)={\rm diag}\left((I_m+X^\tau L)^{-1}X^\tau, \,0\right) \ \ \mbox{for all $X\in {\mathbb{D}}^{m\times n}$.}$$
Therefore, the original $\theta$ is of the form (\ref{35hfhbcnnv7}).

Note that every  surjective homomorphism  from $\mathbb{D}$ to $\mathbb{D}'$ is a ring isomorphism.
Suppose that $\tau$ is surjective.  Then we have $L=0$.  Otherwise, if $L\neq 0$, then there is  $X_0\in{\mathbb{D}}^{m\times n}$
such that $I_m+X_0^\tau L$ is not invertible, a contradiction. Thus, if $\tau$ is surjective, then  $L=0$ and   $\tau$ is a ring isomorphism.
$\qed$

Similarly, we can prove the following proposition by  Lemmas \ref{sd5523w564} and \ref{sdfr533sdv8}.

\begin{proposition}\label{MainLemma002}  Let $\mathbb{D}, \mathbb{D}'$ be division rings with $|\mathbb{D}|\geq 4$,  and let $m,n,m',n'\geq 2$ be
integers with $m'\geq n$, $n'\geq m$. Suppose that  $\theta:  {\mathbb{D}}^{m\times n}\rightarrow  {\mathbb{D}'}^{m'\times n'}$
is a non-degenerate graph homomorphism with  $\theta(0)=0$. Assume further that $\theta(\mathcal{M}_i)\subseteq \mathcal{N}'_i$, $i=1,\ldots,m$, and
$\theta(\mathcal{N}_j)\subseteq \mathcal{M}'_j$, $j=1,\ldots,n$.
Then there exist two invertible diagonal matrices  $P\in GL_{n'}(\mathbb{D}')$ and  $Q\in GL_{m'}(\mathbb{D}')$,
a nonzero ring anti-homomorphism $\sigma: \mathbb{D}\rightarrow\mathbb{D}'$, and a matrix $L\in {\mathbb{D}'}^{m\times n}$
with the property that $I_m+L\,^tX^\sigma\in GL_m(\mathbb{D}')$ for every $X\in {\mathbb{D}}^{m\times n}$, such that
\begin{equation}\label{35hfhbcnnv7002}
 \theta(X)=P\left(
             \begin{array}{cc}
               ^tX^\sigma(I_m+L\,^tX^\sigma)^{-1} & 0 \\
               0 & 0 \\
             \end{array}
           \right)Q, \ \, X\in {\mathbb{D}}^{m\times n}.
\end{equation}
In particular, if $\sigma$ is surjective,  then $L=0$ and   $\sigma$  is a ring anti-isomorphism.
\end{proposition}

Now, we prove Theorem \ref{MainTheorem0} as follows.

\noindent{\bf Proof of Theorem \ref{MainTheorem0}.}  Suppose that $\varphi: {\mathbb{D}}^{m\times n}\rightarrow  {\mathbb{D}'}^{m'\times n'}$
is a non-degenerate graph homomorphism with $\varphi(0)=0$.
Assume further that ${\rm dim}(\varphi(\mathcal{M}_1))=n$ and  ${\rm dim}(\varphi(\mathcal{N}_1))=m$.
We distinguish the following two cases to prove this theorem.

{\bf Case 1.} \ $\varphi(\mathcal{M}_1)$ is contained in a maximal set of type one.
By Lemma \ref{Rectangular-PID2-4} and Lemma \ref{non-degeneratelemma00a}(b),
 there are $P_0\in GL_{m'}(\mathbb{D}')$ and $Q_0\in GL_{n'}(\mathbb{D}')$
such that $\varphi(\mathcal{M}_1)\subseteq P_0{\cal M}_1'Q_0$ and $\varphi(\mathcal{N}_1)\subseteq P_0{\cal N}_1'Q_0$.
Replacing $\varphi$ by the map $X\mapsto P_0^{-1}\varphi(X)Q_0^{-1}$, we get
\begin{equation}\label{vcbsf005hfh}
\varphi({\cal M}_1)\subseteq {\cal M}'_1 \ \ {\rm and} \ \ \varphi({\cal N}_1)\subseteq {\cal N}'_1.
\end{equation}
Let $E_{ij}=E_{ij}^{m\times n}$ and $E_{ij}'=E_{ij}^{m'\times n'}$. It follows from (\ref{vcbsf005hfh}) that
\begin{equation}\label{vcda00oop8}
\mbox{$\varphi(E_{11})=d_{11}E_{11}'$ where $d_{11}\in {\mathbb{D}'}^*$.}
\end{equation}

By Lemma \ref{degenerate-2},  the restriction map $\varphi\mid_{\mathcal{M}_1}: \mathcal{M}_1\rightarrow \mathcal{M}'_1$
is an injective  weighted semi-affine  map.  By  Lemma \ref{weaksemilinea34ear2},
we have
\begin{equation}\label{vxzz1hhp7700}
\varphi\left(
         \begin{array}{c}
           x_1 \\
           0_{m-1,n} \\
         \end{array}
       \right)=
\left(\begin{array}{c}
           (k_1(x_1))^{-1}x_1^{\delta}Q_1 \\
           0_{m'-1,n'} \\
         \end{array}
       \right),  \ x_1\in {\mathbb{D}}^n,
\end{equation}
where  $\delta: \mathbb{D}\rightarrow \mathbb{D}'$ is a nonzero ring homomorphism, $Q_1\in {\mathbb{D}'}^{n\times n'}$,  and $k_1(x_1)$ is of the form
$k_1(x_1)=\sum_{j=1}^nx_{1j}^{\delta}a_{j1}+b_1\neq 0$ for all $x_1=(x_{11}, \ldots, x_{1n})\in {\mathbb{D}}^n$.

By (\ref{vcda00oop8}) and (\ref{vxzz1hhp7700}), it is easy to see that $Q_1$ is of the form $\scriptsize Q_1=\left(
\begin{array}{cc}
q_{11} & 0\\
* & * \\
\end{array}
\right)$ where $q_{11}\neq 0$.
Since $n={\rm dim}(\varphi(\mathcal{M}_1))\leq {\rm rank}(Q_1)$, we get that $n'\geq n$ and ${\rm rank}(Q_1)=n$. There is
$\scriptsize Q_2=\left(
\begin{array}{cc}
q_{11} & 0\\
* & * \\
\end{array}
\right)\in GL_{n'}(\mathbb{D}')$ such that $Q_1=(I_n, 0)Q_2$. Therefore, (\ref{vxzz1hhp7700}) becomes
\begin{equation}\label{vxzz1hhp7700-1}
\varphi\left(
         \begin{array}{c}
           x_1 \\
           0_{m-1,n} \\
         \end{array}
       \right)=
\left(\begin{array}{cc}
           (k_1(x_1))^{-1}x_1^{\delta} &0 \\
           0_{m'-1,n}&0 \\
         \end{array}
       \right)Q_2,  \ x_1\in {\mathbb{D}}^n.
\end{equation}

By Lemma \ref{degenerate-2} again, the restriction map $\varphi\mid_{{\cal N}_1}: {\cal N}_1\rightarrow {\cal N}'_1$
is an injective  weighted semi-affine  map. By  Lemma \ref{weaksemilinea34ear2}, we have
\begin{equation}\label{xcwe0jgool500}
\varphi(y_1, 0,\ldots,0)=\left(P_1y_1^{\mu}(k_1'(y_1))^{-1}, 0, \ldots,0 \right),  \  y_1\in \, ^m\mathbb{D},
\end{equation}
where $\mu: \mathbb{D}\rightarrow \mathbb{D}'$ is a nonzero ring homomorphism, $P_1\in {\mathbb{D}'}^{m'\times m}$, and
$k_1'(y_1)$ is of the form $k_1'(y_1)=\sum_{i=1}^mc_{i1}x_{i1}^{\mu}+d_1\neq 0$ \ for all $y_1= \,^t(x_{11}, \ldots, x_{m1})\in \, ^m\mathbb{D}$.

By (\ref{vcda00oop8}) and (\ref{xcwe0jgool500}), it is clear that $P_1$ is of the form $\scriptsize P_1=\left(
\begin{array}{cc}
p_{11} & *\\
0 & * \\
\end{array}
\right)$ where $p_{11}\neq 0$. Since $m={\rm dim}(\varphi(\mathcal{N}_1))\leq{\rm rank}(P_1)$, we obtain that
$m'\geq m$ and ${\rm rank}(P_1)=m$. There is
$\scriptsize P_2=\left(
\begin{array}{cc}
p_{11} & *\\
0 & * \\
\end{array}
\right)\in GL_{m'}(\mathbb{D}')$ such that $\scriptsize P_1=P_2\left(\begin{array}{c}
I_m\\
 0 \\
 \end{array}
 \right)$. Thus (\ref{xcwe0jgool500}) can be written as
\begin{equation}\label{xcwe0jgool500B}
\varphi(y_1, 0,\ldots,0)=P_2\left(
\begin{array}{cc}
y_1^{\mu}(k_1'(y_1))^{-1} & 0\\
0 & 0_{m'-m, n'-1}\\
\end{array}
\right),  \  y_1\in \, ^m\mathbb{D},
\end{equation}
Note that $p_{11}(k_1(x_1))^{-1}=\left(k_1(x_1)p_{11}^{-1}\right)^{-1}$ and $(k_1'(y_1))^{-1}q_{11}=\left(q_{11}^{-1}k_1'(y_1)\right)^{-1}$.
Without loss of generality, (\ref{vxzz1hhp7700-1}) and (\ref{xcwe0jgool500B}) can be written as
$$
\varphi\left(
         \begin{array}{c}
           x_1 \\
           0_{m-1,n} \\
         \end{array}
       \right)=
P_2\left(\begin{array}{cc}
           (k_1(x_1))^{-1}x_1^{\delta} &0 \\
           0_{m'-1,n}&0 \\
         \end{array}
       \right)Q_2,  \ x_1\in {\mathbb{D}}^n;
$$
$$\varphi(y_1, 0,\ldots,0)=P_2\left(
\begin{array}{cc}
y_1^{\mu}(k_1'(y_1))^{-1} & 0\\
0 & 0_{m'-m, n'-1}\\
\end{array}
\right)Q_2,  \  y_1\in \, ^m\mathbb{D}.$$
Replacing  $\varphi$ by the map $X\mapsto P_2^{-1}\varphi(X)Q_2^{-1}$, we obtain that
\begin{equation}\label{hgkas5FSF3v-00}
\varphi\left(
         \begin{array}{c}
           x_1 \\
           0_{m-1,n} \\
         \end{array}
       \right)=
\left(\begin{array}{cc}
           (k_1(x_1))^{-1}x_1^{\delta} &0 \\
           0_{m'-1,n}&0 \\
         \end{array}
       \right),  \ x_1\in {\mathbb{D}}^n;
\end{equation}
\begin{equation}\label{1phflh64oupv45-00}
\varphi(y_1, 0,\ldots,0)=\left(
\begin{array}{cc}
y_1^{\mu}(k_1'(y_1))^{-1} & 0\\
0 & 0_{m'-m, n'-1}\\
\end{array}
\right),  \  y_1\in \, ^m\mathbb{D}.
\end{equation}

By Corollary \ref{Rectangular-PID2-7},
there are exactly two maximal sets containing $0$ and $E_{i1}$  [resp. $0$ and $bE_{i1}'$ where $b\in {\mathbb{D}'}^*$] in ${\mathbb{D}}^{m\times n}$ [resp. ${\mathbb{D}'}^{m'\times n'}$],
they are $\mathcal{N}_1$ and $\mathcal{M}_i$ [resp. $\mathcal{N}'_1$ and $\mathcal{M}'_i$]. Applying  Lemma \ref{non-degeneratelemma00a}(b) and
(\ref{1phflh64oupv45-00}),  we have that $\varphi(\mathcal{N}_1)\subseteq \mathcal{N}_1'$ and
$$\varphi(\mathcal{M}_i)\subseteq \mathcal{M}_i', \ i=1,\ldots,m.$$
Similarly, from (\ref{hgkas5FSF3v-00}) we have
$$\varphi(\mathcal{N}_j)\subseteq \mathcal{N}_j', \ j=1,\ldots,n.$$

By Proposition \ref{MainLemma001}, there exist two invertible diagonal matrices  $T_1\in GL_{m'}(\mathbb{D}')$ and  $T_2\in GL_{n'}(\mathbb{D}')$,
a nonzero ring homomorphism $\tau: \mathbb{D} \rightarrow\mathbb{D}'$, and a matrix $L\in {\mathbb{D}'}^{n\times m}$ with the property
that $I_m+X^\tau L\in GL_m(\mathbb{D}')$ for every $X\in {\mathbb{D}}^{m\times n}$, such that
$$
\varphi(X)=T_1\left(\begin{array}{cc}
 (I_m+X^\tau L)^{-1}X^\tau  & 0 \\
  0 & 0 \\
   \end{array}
  \right)T_2, \ \, X\in {\mathbb{D}}^{m\times n}.$$
Therefore, the original $\varphi$ is of the form (\ref{tt3654665cnc00}).
In particular, if $\tau$ is surjective,  then $L=0$ and  $\tau$ is a ring isomorphism.

{\bf Case 2.} \ $\varphi(\mathcal{M}_1)$ is contained in a maximal set of type two.
Similar to the Case 1, replacing the map $\varphi$ by one map of the form $X\mapsto {P'}^{-1}\varphi(X){Q'}^{-1}$,
where $P',Q'$ are two invertible matrices over $\mathbb{D}'$,
we have that  $m'\geq n$, $n'\geq m$, $\varphi(\mathcal{M}_i)\subseteq \mathcal{N}_i'$, $i=1,\ldots,m$, and
$\varphi(\mathcal{N}_j)\subseteq \mathcal{M}_j'$, $j=1,\ldots,n$.
By Proposition \ref{MainLemma002}, the original $\varphi$ is of the form (\ref{ttd2CX3mmbb6600}).
In particular, if $\sigma$ is surjective,  then $L=0$ and  $\sigma$ is a ring anti-isomorphism.  \hfill $\Box$

Finally, we prove Corollary \ref{MainTheoremcor0} as follows.

\noindent{\bf Proof of Corollary \ref{MainTheoremcor0}.} \
{\em Case 1}. \ $\varphi(\mathcal{M}_1)$ is contained in a maximal set of type one. Then there is  $P'\in GL_{m'}(\mathbb{D})$
such that  $\varphi(\mathcal{M}_1)\subseteq P'\mathcal{M}_1'$.
 By Lemmas  \ref{degenerate-2} and \ref{weaksemilinea34ear2}, the restriction map $\varphi\mid_{\mathcal{M}_1}: \mathcal{M}_1\rightarrow P'\mathcal{M}_1'$ is of the form
$$
\varphi\left(
         \begin{array}{c}
           x \\
           0 \\
         \end{array}
       \right)=
P'\left(\begin{array}{c}
(k(x))^{-1}x^{\delta}T \\
0 \\
\end{array}
\right),  \  x\in {\mathbb{D}}^n,
$$
where  $\delta: \mathbb{D}\rightarrow \mathbb{D}$ is a nonzero endomorphism, $T\in {\mathbb{D}}^{n\times n'}$, and
$k(x)=\sum_{j=1}^nx^{\delta}_{1j} a_{j1}+b\neq 0$ for all $x=(x_{11}, \ldots, x_{1n})\in {\mathbb{D}}^n$.
Since $\mathbb{D}$ is an EAS division ring, $\delta$ is bijection and hence it is a ring automorphism of $\mathbb{D}$. Thus, we must have ${\rm rank}(T)=n$.
Otherwise, there exists  $0\neq \alpha\in{\mathbb{D}}^n$ such that $(k(\alpha))^{-1}\alpha^{\delta}T=0$, a contradiction to
$\varphi\mid_{\mathcal{M}_1}$ being injective.
Therefore $n'\geq n$.

Let $e_i$ and $\beta_i$ be the $i$-th rows of $I_n$ and $T$, respectively.
  Since
$\scriptsize\varphi\left(
         \begin{array}{c}
           e_{i} \\
           0 \\
         \end{array}
       \right)=
P'\left(\begin{array}{c}
(k(e_{i}))^{-1}\beta_{i} \\
0 \\
\end{array}
\right)$ and ${\rm rank}(T)=n$,  matrices
$\scriptsize{P'}^{-1}\varphi\left(
         \begin{array}{c}
           e_{1} \\
           0 \\
         \end{array}
       \right), \ldots, {P'}^{-1}\varphi\left(
         \begin{array}{c}
           e_{n} \\
           0 \\
         \end{array}
       \right)$ are $n$
 left linearly independent matrices over $\mathbb{D}$, which implies that ${\rm dim}(\varphi(\mathcal{M}_1))\geq n$.
By Lemma \ref{degenerate-2}, ${\rm dim}(\varphi(\mathcal{M}_1))\leq{\rm dim}(\mathcal{M}_1)=n$. Therefore,
${\rm dim}(\varphi(\mathcal{M}_1))=n$.

By Lemma \ref{non-degeneratelemma00a}(b) and Lemma \ref{Rectangular-PID2-4},  there exists  $Q'\in GL_{n'}(\mathbb{D})$
such that $\varphi(\mathcal{N}_1)\subseteq \mathcal{N}_1'Q'$. Similarly, we can prove that $m'\geq m$ and
${\rm dim}(\varphi(\mathcal{N}_1))=m$.
By Theorem \ref{MainTheorem0}, there exist  matrices  $P\in GL_{m'}(\mathbb{D})$ and $Q\in GL_{n'}(\mathbb{D})$, and
a ring automorphism $\tau$ of $\mathbb{D}$,  such that
$$\varphi(X)=P\left(\begin{array}{cc}
X^\tau  & 0 \\
  0 & 0 \\
   \end{array}
  \right)Q, \ \, X\in {\mathbb{D}}^{m\times n}.$$

{\em Case 2}. \ $\varphi(\mathcal{M}_1)$ is contained in a maximal set of type two.  Then there is $Q'\in GL_{n'}(\mathbb{D})$
such that $\varphi(\mathcal{M}_1)\subseteq \mathcal{N}_1'Q'$.
Similarly, we can prove that $m'\geq n$, $n'\geq m$, ${\rm dim}(\varphi(\mathcal{M}_1))=n$
 and ${\rm dim}(\varphi(\mathcal{N}_1))=m$.
By Theorem \ref{MainTheorem0}, there exist  matrices  $P\in GL_{m'}(\mathbb{D}')$ and $Q\in GL_{n'}(\mathbb{D}')$, and
a ring anti-automorphism $\sigma$ of $\mathbb{D}$,  such that
$\scriptsize\varphi(X)=P\left(\begin{array}{cc}
^tX^\sigma  & 0 \\
  0 & 0 \\
   \end{array}
  \right)Q, \ \, X\in {\mathbb{D}}^{m\times n}$. \hfill $\Box$



\end{document}